\documentclass{tran-l}

\usepackage{amsmath}
\usepackage{amsfonts}
\usepackage{amssymb}

\newcommand{\abs}[1]{\ensuremath{\left| #1 \right| }}
\newcommand{\qPs}{$q$--Pochhammer symbol}
\newcommand{\bhs}{basic hypergeometric series}
\newcommand{\rhs}{right hand side}
\newcommand{\lhs}{left hand side}
\newcommand{\wps}{well--poised side}
\newcommand{\bs}{balanced side}
\newcommand{\wrt}{with respect to}
\newcommand{\od}{one dimensional}
\newcommand{\hd}{higher dimensional}
\newcommand{\rc}{rectangular case}
\newcommand{\BT}{Bailey Transform}
\newcommand{\BL}{Bailey Lemma}
\newcommand{\Bl}{Bailey lattice}
\newcommand{\Bc}{Bailey chain}
\newcommand{\Bp}{Bailey pair}
\newcommand{\Wt}{Watson transformation}
\newcommand{\RRi}{Rogers--Ramanujan identity}
\newcommand{\RRis}{Rogers--Ramanujan identities}
\newcommand{\AGis}{Andrews--Gordon identities}
\newcommand{\Epnt}{Euler's Pentagonal Number Theorem}
\newcommand{\RSi}{Rogers--Selberg identity}
\newcommand{\RSis}{Rogers--Selberg identities}
\newcommand{\Mis}{Macdonald identities}
\newcommand{\Wg}{Weyl group}
\newcommand{\Wdegf}{Weyl degree formula}
\newcommand{\Wdenf}{Weyl denominator formula}
\newcommand{\hgr}{hyperoctahedral group}
\newcommand{\Js}{Jackson sum}
\newcommand{\ci}{cocycle identity}
\newcommand{\Jtpi}{Jacobi triple product identity} 
\newcommand{\tns}{\ensuremath{_{10}\varphi_9}}
\newcommand{\sfs}{\ensuremath{_6\varphi_5}}
\newcommand{\tts}{\ensuremath{_3\varphi_2}}
\newcommand{\ml}{multilateralization lemma}
\newcommand{\Jc}{Jackson coefficients}
\newcommand{\eJc}{elliptic Jackson coefficients}

\newtheorem{theorem}{Theorem}[section]
\newtheorem{lemma}[theorem]{Lemma}
\newtheorem{cor}[theorem]{Corollary} 

\theoremstyle{definition}
\newtheorem{definition}[theorem]{Definition}

\theoremstyle{remark}
\newtheorem{remark}[theorem]{Remark}

\numberwithin{equation}{section}

\begin{document}

\title[Multiple $q$--Series Identities] 
{An Elliptic $BC_n$ Bailey Lemma, \\ 
Multiple Rogers--Ramanujan Identities and \\
Euler's Pentagonal Number Theorems}   

\author{Hasan Coskun}  
\address{Department of Mathematics, Texas A\&M
  University--Commerce, Binnion Hall, Room 314, Commerce, TX 75429}  
\email{hasan\_coskun@tamu-commerce.edu}

\date{May 21, 2006 and, in revised form, October 14, 2006.}

\begin{abstract}
An elliptic $BC_n$ generalization of the classical two parameter
Bailey Lemma is proved, and a basic one parameter $BC_n$ Bailey Lemma
is obtained as a limiting case. Several summation and transformation
formulas associated with the root system $BC_n$ are proved as
applications, including a $_6\varphi_5$ summation formula, a
generalized Watson transformation and an unspecialized Rogers--Selberg
identity. The last identity is specialized to give an infinite family
of multilateral Rogers--Selberg identities. Standard determinant
evaluations are then used to compute $B_n$ and $D_n$ generalizations
of the Rogers--Ramanujan identities in terms of determinants of theta
functions. Starting with the $BC_n$ $_6\varphi_5$ summation formula, a
similar program is followed to prove an infinite family of $D_n$
Euler's Pentagonal Number Theorems.   
\end{abstract}



\maketitle

\section{Introduction}
\label{section1}
The \RRis\ and \Epnt\ are decisively among the most celebrated
classical $q$--series identities. These identities are usually written
in terms of the \qPs\ $(a;q)_\alpha$ for $q, \alpha\in\mathbb{C}$,
which is defined formally as  
\begin{equation}
\label{qPochSymbol}
(a)_\alpha = (a;q)_\alpha :=\dfrac{(a;q)_\infty}{(aq^\alpha;q)_\infty}
\end{equation} 
where $(a;q)_\infty:=\prod_{i=0}^{\infty} (1-aq^i)$. 

With this notation, the \RRis\ can be written in the form
\begin{equation}
  \label{eq:RRs}
  \sum_{m=0}^\infty \dfrac{q^{m(m+\delta)}}{(q; q)_m} =
  \dfrac{1}{(q^{1+\delta}; q^5)_\infty (q^{4-\delta}; q^5)_\infty}
\end{equation}
where $\delta\in\{0,1\}$ and $\abs{q}<1$. These identities have a very
rich history. Many important figures in mathematics had contributed to
the development of these identities starting with
Rogers~\cite{Rogers1} who first proved them in 1894, and
Ramanujan~\cite{Hardy1} whose involvement made Rogers' unnoticed work
popular. Others contributed by simplifying existing proofs, suggesting
new proofs of different nature, establishing their relations to other
branches of mathematics and generalizing these
identities~\cite{AndrewsS1}, \cite{Bailey1}, \cite{GasperR1},
\cite{Milne7}, \cite{Slater2}, \cite{Stembridge1}.  This paper proves
multiple series analogues of the \RRis\ associated to root systems
$B_n$ and $D_n$ of rank~$n$.

Section \S\ref{section5} reviews definitions of two remarkable
families of $BC_n$ type symmetric rational functions~\cite{Coskun},
namely Macdonald functions $W_{\lambda/\mu}(z;r,q,p,t; a,b)$ and \Jc\
$\omega_{\lambda/\mu}(z; r,q,p,t; a,b)$. Properties of these
functions are then used to prove an elliptic $BC_n$ generalization of
the classical two parameter \BL. In particular, an important elegant
property called the \ci\ for $\omega_{\lambda/\mu}$ played an
important role in this program. It may be written in the form   
\begin{multline}
\label{cocycleIdentity}
  \omega_{\nu/\mu}((uv)^{-1};uv,q,p,t;a(uv)^2, buv) \\ = \sum_{\mu\subseteq
  \lambda \subseteq \nu} \omega_{\nu/\lambda}(v^{-1};v,q,p,t;a(vu)^2,bvu) \,
  \omega_{\lambda/\mu}(u^{-1};u,q,p,t;au^2,bu)   
\end{multline}
where the summation index $\lambda$ runs over partitions. 

A basic one parameter $BC_n$ \BL\ is obtained as a limiting case of
the two parameter $BC_n$ \BL\ in the same section. The one--parameter \BL\
is iterated to generate several remarkable $BC_n$ analogues of the
classical \bhs\ identities including a
$_6\varphi_5$ summation formula and a generalized Watson transformation.

Section \S\ref{section4} gives a $BC_n$ \RSi\ as a limiting case of
the $BC_n$ \Wt. Specializing the parameters of this identity gives
rise to an infinite family of multilateral Rogers--Selberg
identities associated with the root systems $B_n$ and $D_n$, using a very
general multilateralization argument. Standard 
determinant evaluations are then used to compute the
$B_n$ and $D_n$ generalizations of the Rogers--Ramanujan identities in 
terms of determinants of theta functions.

With the notation as above, $D_n$ multiple \RRis\ can be written as
\begin{multline}
\label{intro:DnRRis}
\sum_{\lambda \in \mathbb{Z}_{\geq}^n } \prod_{i=1}^{n} \left\{ \dfrac{
q^{(\delta+n-1) (\lambda_i-n+i) + (\lambda_i-n+i)^2 } } {(q
)_{\lambda_i }} \right\} \prod_{1\leq i<j \leq n} \left\{
(1-q^{\lambda_i-\lambda_j} ) \right\} \\ = \dfrac{(-1)^{\binom{n}{2}}
q^{\binom{n}{2} (4n+3\delta-2)/6} }{2 \, \theta(q;q^5)^n \,
\theta(q^2;q^5)^n} \cdot \!\!\! \det_{1\leq i,j\leq n} \bigg(
q^{(j-1)(n-i+\delta/2)}\, \theta(q^{4n+2\delta+1-4i+j}; q^5 ) \bigg.
\bigg. \\ + q^{-(j - 1) (n-i+\delta/2)} \, \theta(q^{4n+2\delta
+3-4i-j} ;q^5) \bigg)
\end{multline}
where $\theta(z;q):= (z;q)_\infty (q/z;q)_\infty$. The cases $\delta=0$ and
$\delta=1$ give the first and the second $D_n$ \RRis, respectively.  A
single $B_n$ multiple \RRi\ can be written, similarly, in the form
\begin{multline}
\label{intro:BnRRis}
\sum_{\lambda \in \mathbb{Z}_{\geq}^n } \prod_{i=1}^{n} \left\{ \dfrac{
q^{(1+n) (\lambda_i-n+i) + (\lambda_i-n+i)^2 } } {(q )_{\lambda_i}}
\right\} \prod_{1\leq i<j \leq n} \left\{ (1-q^{\lambda_i-\lambda_j} )
\right\} \\ = \dfrac{(-1)^{\binom{n}{2}+n} q^{n(n+1)(4n-1)/12}
}{\theta(q;q^5)^n \, 
\theta(q^2;q^5)^n } \cdot \det_{1\leq i,j\leq n} \bigg(
q^{(j-1/2)(n-i+1)} \theta(q^{6+4n-4i+j}; q^5 ) \\ - q^{-(j -
1/2)(n-i+1)} \theta(q^{7+4n-4i-j}; q^5 ) \bigg)
\end{multline}
In both cases, $n$ is a positive integer and $\,\abs{q}<1$ as
usual. Alternative versions of~(\ref{intro:DnRRis})
and~(\ref{intro:BnRRis}) in form of determinant transformation
identities are given in~(\ref{RRDetTransfIdenDn})
and~(\ref{RRDetTransfIdenBn}), respectively.

The classical \Epnt\ states that 
\begin{equation}
  \label{eq:Epnt}
  (q)_\infty = \sum_{m=0}^\infty (-1)^m q^{m+3\binom{m}{2}} 
\end{equation}
This beautiful identity was first proved by Euler~\cite{Euler}
in 1793. Many generalizations and combinatorial interpretations 
appeared in literature~\cite{Andrews3}, \cite{Andrews4},
\cite{BerkovichG1}, \cite{Dyson}, \cite{SchillingW1}, etc.  This paper
gives a remarkable multiple series 
analogue of \Epnt\ associated to the root system $D_n$ of rank~$n$. In
fact, Section \S\ref{section4} proves an infinite family of \Epnt s
using the $BC_n$ \sfs\ summation formula from Section
\S\ref{section5}, and methods developed in that section. These
identities can be written in the form
\begin{multline}
\label{intro:Epnt}
(q; q)^n_\infty \!\!\! \prod_{1\leq i<j\leq n} \dfrac{(q^{k(j-i)}; q)_\infty
  }{(q^{k(j-i+1)} ; q)_{\infty} } = \sum_{\mu\in\mathbb{Z}^n }
  (-1)^{|\mu|} q^{-kn(\mu)+3n(\mu')+|\mu|(k(n-1)+1)} \\ \cdot
  \sideset{}{^{k,1}}\prod_{1\leq i<j \leq n} \left\{
  \dfrac{(q^{1+k(j-i-1) + \mu_i-\mu_j}; q)_{\infty} } {(q^{k(1+2n-i-j) +
  \mu_i+\mu_j}; q)_{\infty} } \dfrac{ (q^{1+k(-1+2n-i-j) +
  \mu_i+\mu_j}; q)_{\infty} } {(q^{k(j-i+1) + \mu_i -\mu_j};
  q)_{\infty} } \right\}
\end{multline}
where $n\in\mathbb{Z}_{>}$, $m, k\in\mathbb{Z}_{\geq}$, and 
$\,\abs{q}<1$ as usual.
Similar to~(\ref{RRDetTransfIdenDn}) and~(\ref{RRDetTransfIdenBn}),
the $k=1$ instance of~(\ref{intro:Epnt}) is written as a determinant
identity in terms of theta functions in~(\ref{Epnt_k1}).

Some results of this paper, such as~(\ref{intro:DnRRis}), were
conjectured in author's PhD work~\cite{Coskun} which 
was conducted under the supervision of Dr. R. A. Gustafson. 

Section~\ref{AGeneralRRis} outlines an alternative proof of a known
result. This proof may be seen as a one dimensional
version of the more general proofs used in higher
dimensional results in later sections.  It is hoped that this
section helps to make the rest of the paper more readable.

\section{A Generalization of \RRis}
\label{AGeneralRRis}
Garrett, Ismail and Stanton~\cite{GarretI1} used orthogonal polynomials to
obtain, among other results, the following generalization of the \RRis. 
\begin{multline}
\label{eq:GISgen}
\sum_{\lambda_i=0}^\infty
  \dfrac{q^{\lambda_i(\lambda_i+\delta)}}{(q; q)_{\lambda_i}} =
  \dfrac{(-1)^\delta q^{-\binom{\delta}{2} } E_{\delta-2}(q) } {(q,
  q^{4}; q^5)_\infty} - \dfrac{(-1)^{\delta} q^{-\binom{\delta}{2} }
  D_{\delta-2}(q) } {(q^2, q^{3}; q^5)_\infty} \\
= (-1)^\delta q^{-\binom{\delta}{2} } \dfrac{(q^5; q^5)_\infty }
  {(q)_\infty } \left[ \theta(q^2; q^5) E_{\delta-2}(q) - \theta(q; q^5) 
  D_{\delta-2}(q)  \right] 
\end{multline}
where $\delta\in\mathbb{Z}$ (also see~\cite{BerkovichP1}), and the
Schur polynomials $D_\delta(q)$ and $E_\delta(q)$ are defined by
\begin{eqnarray}
\label{eq:SchurPolyRecursive}
D_\delta(q) &=& D_{\delta-1} + q^\delta D_{\delta-2}, \quad D_0=1, D_1=1+q \\
E_\delta(q) &=& E_{\delta-1} + q^\delta E_{\delta-2}, \quad E_0=1, E_1=1 
\end{eqnarray}
recursively. In two special cases $\delta\in\{1,0\}$ the
formula~(\ref{eq:GISgen}) simplifies to give the classical
\RRis~(\ref{eq:RRs}). 

There are numerous (analytic, combinatorial, probabilistic, algebraic)
proofs~\cite{Andrews2} of the
\RRis~(\ref{eq:RRs}). Watson's proof~\cite{Watson1} of 
these identities, for example, depends on the identity 
\begin{multline}
\label{eq:onedimWt}
_8\varphi_7 \left[
\begin{matrix} 
b, q b^{1/2}, -q b^{1/2}, \sigma_1, \rho_1 ,\sigma_2, \rho_2, q^{-N} \\
b^{1/2}, -b^{1/2}, bq/\sigma_1, bq/\rho_1, bq/\sigma_2, bq/\rho_2,
bq^{N+1} 
\end{matrix}
;q , \dfrac{b^2q^{2+N}}{\sigma_1\sigma_2\rho_1\rho_2} \right] \\ =
 \dfrac{(bq, bq/\sigma_2 \rho_2)_N} {(bq/\sigma_2, 
 bq/\rho_2)_N } {_4\varphi_3} \left[
\begin{matrix} 
bq/\sigma_1 \rho_1, \sigma_2, \rho_2 , q^{-N} \\
bq/\sigma_1,  bq/\rho_1, \sigma_2 \rho_2 q^{-N}/b 
\end{matrix}
;q , q \right]
\end{multline}
called \Wt. Watson showed that in the limiting case as $N$, $\sigma_1$,
$\sigma_2, \rho_1,\rho_2 \rightarrow \infty$, the
transformation~(\ref{eq:onedimWt}) gives the remarkable \RSi   
\begin{equation}
\label{eq:onedimRSi}
\sum_{m=0}^\infty \dfrac{b^m q^{m^2} } {(q)_m} =
\dfrac{1}{(qb)_\infty } \sum_{m=0}^\infty  (-1)^m b^{2m }
q^{m(5m-1)/2} \dfrac{(1-bq^{2m} )}{(1-b) }
\dfrac{(b)_m }{(q)_m } 
\end{equation} 
In special cases $b=q^\delta$ for $\delta\in\{1,0\}$, the
series on the \rhs\ of~(\ref{eq:onedimRSi}) can be written as a
bilateral sum giving  
\begin{equation}
\label{eq:onedimRSiBilateral}
  \sum_{m=0}^\infty \dfrac{q^{m(m+\delta)} } {(q)_m} =
 \dfrac{1}{(q)_\infty } \sum_{m=-\infty}^\infty (-1)^m
 q^{5\binom{m}{2}+2(\delta+1)m} 
\end{equation} 
The product representation~(\ref{eq:RRs}) can now be computed applying
the \Jtpi\ 
\begin{equation}
\label{eq:JtpiBackground}
\theta(z;q) := (z,q/z;q)_\infty = \dfrac{1}{(q;q)_\infty }
\sum_{m=-\infty}^\infty (-1)^{m} q^{\binom{m}{2}} z^m, \quad
\abs{q}<1
\end{equation}
to the \rhs\ of~(\ref{eq:onedimRSiBilateral}) after rescaling parameters
$q$ by $q^5$ and $z$ by $q^{2(\delta+1)}$ for $\delta\in\{1,0\}$. 

In this section an alternative elementary proof of the
identity~(\ref{eq:GISgen}) will be given using the symmetries of
the \RSi~(\ref{eq:onedimRSi}).
A generalization of this argument will then be used in the proof of
our multiple \RRis\ below.   

By setting $b=q^{2z}$ for some $z\in\mathbb{C}$, flipping
appropriate factors using the definition of the
\qPs~(\ref{qPochSymbol}) and simplifying,  
the \od\ \RSi~(\ref{eq:onedimRSi}) can be written in the form 
\begin{equation}
\label{eq:onedimRSiBilateral2V}
\sum_{m=0}^\infty q^{(z+m)^2} (q^{1-z+(z+m)})_{\infty}  
= \sum_{m=0}^\infty q^{(z+m)^2} \dfrac {(q^{1-z+(z+m)})_{\infty}
  (q^{1-z-(z+m)})_{\infty }} {(q^{1+2(z+m)})_\infty  
  (q^{1-2(z+m)})_\infty }
\end{equation}
The summand on the \wps\ is obviously invariant under the action of
the \Wg\ for $C_1$, that is $W=\mathbb{Z}_2=\{\imath, w\}$, generated by
the maps  
$q^{z+m} \leftrightarrow q^{\pm(z+m)}$ for any $z\in\mathbb{C}$. This
implies that 
\begin{equation}
\sum_{m\in L_{+} } f(m,z)  = \sum_{m\in wL_{+}} f(m,z)
\end{equation}
where $f$ is the summand in the \rhs\ of~(\ref{eq:onedimRSiBilateral2V}),
$L_+ = \imath L_+$ is the lattice of all non--negative integers and 
$wL_{+} = \{ -m-2z: m\in L_+\} $. 

The \rhs\ of the classical \RSi~(\ref{eq:onedimRSi}) can be written in an
equivalent form by flipping different terms as 
\begin{equation}
\label{eq:onedimRSi2ndversion}
\dfrac{1} {(q)_{\infty }}  \sum_{m=0}^\infty q^{4zm} (-1)^{m}
   q^{2m+5\binom{m}{2}} 
   \dfrac {(q^{2(z+m)})_\infty } {(q^{1+2(z+m)})_\infty} \dfrac
   {(q^{1+m})_{\infty} } { (q^{z+(z+m)})_\infty }   
\end{equation}
Now multiply the series~(\ref{eq:onedimRSi2ndversion}) by the \rhs\
of the 
$C_1$ version of Macdonald's~\cite{Macdonald4} polynomial identity
(written slightly differently as)
\begin{equation}
\label{C1MacdIdent}
 1= \dfrac{(q^{1+2z+2m})_\infty } {(q^{2z+2m})_\infty} +
 \dfrac{(q^{1-2z-2m})_\infty } {(q^{-2z-2m})_\infty} 
\end{equation}
Since terms of the sum~(\ref{eq:onedimRSi2ndversion}) as well as both
terms of~(\ref{C1MacdIdent}) are invariant under the maps $q^{z+m}
\leftrightarrow q^{\pm(z+m)}$, the
series~(\ref{eq:onedimRSi2ndversion}) can be written in the form
\begin{equation}
\label{eq:onedimRSi2ndversionReduced}
\dfrac{1} {(q)_{\infty }}  \sum_{m\in L_z} 
   (-1)^{m} q^{2(1+2z)m+5\binom{m}{2}} 
\dfrac{(q^{1+m})_{\infty} } { (q^{2z+m})_\infty }   
\end{equation}
where $L_z = \cup_{w\in W} wL_+$. It is clear that the
series~(\ref{eq:onedimRSi2ndversionReduced}) is possibly over the full
weight lattice 
$L=\mathbb{Z}$ for $C_1$ (i.e., $L_z=L$) only when $z=\delta/2$ for some
integer $\delta$. In that case, one still needs to study possible
overlaps and gaps between the lattices $L_+$ and $wL_+$. 

Let's denote by $\mathcal{O}$ the ``overlap'' set 
$L_+ \cap wL_+$  and by $\mathcal{G}$ the ``gap'' $L
\backslash L_z$. Clearly, the overlap $\mathcal{O} \neq
\emptyset$ and $\mathcal{G} = \emptyset$ if $\delta$ is non--positive,
and $\mathcal{O} = \emptyset$ and $\mathcal{G} \neq
\emptyset$ otherwise.   

Assume that $\delta>0$. The
series~(\ref{eq:onedimRSi2ndversionReduced}) can 
be written over $L$, since the summand $g(m,z)$
in~(\ref{eq:onedimRSi2ndversionReduced}) vanishes on $\mathcal{G} =
\{-1,\ldots, -2z+1\}=\{-1,\ldots, -\delta+1\}$ due to the numerator factor
$(q^{1+m})_{\infty}$. Note that the denominator factor
$(q^{2z+m})_\infty\neq 0$ on $\mathcal{G}$. In particular, the second
\RRi\ corresponding to $\delta=1$ follows from the case $z=1/2$. 

For $\delta\leq 0$, the additional ``overlap
condition'' needs to be verified. Namely,
\begin{equation}
\label{overlapCondOneDim}
 \sum_{m\in\mathcal{O}} f(m,z) = \sum_{m\in\mathcal{O}} 
g(m,z)
\end{equation}
where $\mathcal{O} = \{0, \ldots, -2z\}= \{0, \ldots, -\delta \}$. In
particular, for  
$z=0$ the overlap is $\mathcal{O} = \{0\}$ and it is plain that
$f(0,0)=g(0,0)$. The first \RRi\ then follows 
using a rewriting of the $BC_1$ version of 
Macdonald's~\cite{Macdonald4} polynomial identity 
\begin{equation}
\label{BC1MacdIdent}
 \dfrac{(-q^{1+z+m})_\infty } {(-q^{z+m})_\infty} +
 \dfrac{(-q^{1-z-m})_\infty } {(-q^{-z-m})_\infty} =1
\end{equation}

Multiplying the specialized \RSi\ by Macdonald's polynomial
identities~(\ref{C1MacdIdent}) and~(\ref{BC1MacdIdent}) amounts to
dropping certain factors corresponding to positive (or
negative) roots in the summand. 

It is already verified that the \RSi~(\ref{eq:onedimRSi})
can be written as a bilateral sum when $b=q^\delta$ for all
$\delta\in\mathbb{Z}_{\geq}$.  
Now defining a degree $\delta+1$ polynomial $f_\delta$ in $x$ by
\begin{equation}
f_\delta(x):=(1-x^2 q^{\delta}) (qx)_{\delta-1} \quad \mathrm{and}
\quad f_0(x):=(1+x),
\end{equation}
and flipping appropriate terms,   
the \rhs\ of the \RSi~(\ref{eq:onedimRSi}) may be written in the form 
\begin{equation}
\dfrac{1}{2 (q)_{\delta}} \sum_{m=-\infty}^\infty (-1)^{m}
q^{2(1+\delta)m+5\binom{m}{2}} f_\delta(q^m)
\end{equation}
Using the \Jtpi~(\ref{eq:JtpiBackground}) one finally gets
\begin{equation}
\sum_{m=0}^\infty \dfrac{ q^{m(m+\delta )} } {(q)_{m}} 
= \dfrac{(q^5; q^5)_\infty}{2 (q;q)_{\infty} } \sum_{n=0}^{\delta+1}
   \dfrac{f_\delta^{(n)}(0)}{n!} \, \theta(q^{2(\delta+1)+n}; q^5)
\end{equation}
The coefficients $f_\delta^{(n)}(0) /n!$ can be easily
computed via the terminating $q$-binomial theorem
\begin{equation}
(x;q)_{\delta}= \sum_{m=0}^\delta \left[ \begin{matrix} \delta \\
   m \end{matrix} \right]_q (-1)^m q^{\binom{m}{2}} x^m 
\end{equation}
in terms of $q$-binomial coefficients defined by 
\begin{equation}
\left[\begin{matrix} n \\ m \end{matrix} \right]_q := \dfrac{(q)_n
}{(q)_{n-m} (q)_m}  
\end{equation}
when $0\leq m\leq n$, and 0 otherwise. The coefficients of $f_\delta$
yields a well--known alternative representation 
of Schur polynomials~(\ref{eq:SchurPolyRecursive}) 
\begin{equation}
E_{\delta-2}(q) = \sum_{k} (-1)^k q^{k (5k-3)/2} 
\left[ \begin{matrix} \delta-1 \\ \lfloor \frac{\delta+1-5k}{2} \rfloor
\end{matrix} \right]_q 
\end{equation}
and
\begin{equation}
D_{\delta-2}(q) = \sum_{k} (-1)^k q^{k (5k+1)/2} 
\left[ \begin{matrix} \delta-1 \\ \lfloor \frac{\delta-1-5k}{2} \rfloor
\end{matrix} \right]_q .
\end{equation}
for $\delta\geq 2$.

A slight generalization of this argument shows that~(\ref{eq:GISgen})
also holds when $\delta<0$, and completes the proof
of~(\ref{eq:GISgen}).

\section{An Elliptic $BC_n$ Bailey Lemma}
\label{section5}
In this section, an elliptic $BC_n$ \BL\ will be proved. This result
will then be used in the next section to obtain multiple \RRis\ associated to
root systems.

First, definitions are given for infinite dimensional matrices $M(a,b)$
and $S(b)$ indexed by
partitions \wrt\ partial inclusion ordering $\subseteq$ defined by
\begin{equation}
\label{partialordering}
\mu \subseteq \lambda \;\Leftrightarrow \;\mu_i \leq \lambda_i, \quad
\forall i\geq 1.
\end{equation}
As in~\cite{Coskun1}, the $\mathbb{Z}$--space $V$ denotes the space of
infinite lower--triangular matrices whose entries are rational
functions in complex parameters
$\rho_i,\sigma_i$ for $i\in\mathbb{Z}_{\geq}$ over the field
$\mathbb{F}=\mathbb{C}(q, p, t,r,a,b)$. 
The condition that $u\in
V$ is lower triangular \wrt\ the partial inclusion
ordering~(\ref{partialordering}) can be stated in the form
\begin{equation}
  u_{\lambda\mu} = 0 ,\, \quad \mathrm{when}\; \mu \not \subseteq
  \lambda.
\end{equation}
The multiplication operation in $V$ is defined by the relation
\begin{equation}
\label{multiplication}
  (uv)_{\lambda\mu} := \sum_{\mu\subseteq\nu\subseteq\lambda}
  u_{\lambda\nu} v_{\nu\mu}
\end{equation}
for $u,v\in V$. 

The definitions of $M(a,b)$ and $S(b)$ involve the symmetric elliptic
Macdonald functions $W_{\lambda/\mu}$ and \Jc\ $\omega_{\lambda/\mu}$
on $BC_n$ defined and investigated in~\cite{Coskun}
and~\cite{Coskun1}. A brief review of the definitions and basic
properties of these functions are in order. 

Recall that an elliptic
analogue of the basic factorial is given in terms of 
$\theta(x)$ function as follows~\cite{FrenkelT1}. For $x, p\in
\mathbb{C}$ and $\abs{p}<1$, let
\begin{equation}
  \theta(x) = \theta(x;p) := (x; p)_\infty (p/x; p)_\infty
\end{equation}
and for $a\in \mathbb{C}$, and a positive integer $m$ define
\begin{equation}
  (a; q,p)_m := \prod_{k=0}^{m-1} \theta(aq^m)
\end{equation}
The definition is extended to negative $m$ by setting
  $(a; q,p)_m = 1/ (aq^{m}; q, p)_{-m} $. 
Note also that when $p=0$, $(a; q,p)_m$ reduces to standard
(trigonometric) \qPs.

For any partition $\lambda = (\lambda_1, \ldots, \lambda_n)$ and
$t\in\mathbb{C}$, one can also define~\cite{Warnaar2}
\begin{equation}
\label{ellipticQtPocSymbol}
  (a)_\lambda=(a; q, p, t)_\lambda := \prod_{k=1}^{n} (at^{1-i};
  q,p)_{\lambda_i} .
\end{equation}
Note that when $\lambda=(\lambda_1) = \lambda_1$ is a single part
partition, then $(a; q, p, t)_\lambda = (a; q, p)_{\lambda_1} =
(a)_{\lambda_1}$. The following notation will also be used.
\begin{equation}
  (a_1, \ldots, a_k)_\lambda = (a_1, \ldots, a_k; q, p, t)_\lambda :=
  (a_1)_\lambda \ldots (a_k)_\lambda .
\end{equation}

Now let $\lambda=(\lambda_1, \ldots, 
\lambda_n)$ and $\mu=(\mu_1, \ldots, \mu_n)$ be partitions of at most
$n$ parts for a positive integer $n$ such that the
skew partition $\lambda/\mu$ is a horizontal strip; i.e. $\lambda_1
\geq \mu_1 \geq\lambda_2 \geq \mu_2 \geq \ldots \geq \lambda_n \geq
\mu_n$. Then for $q,p,t,x,a,b\in\mathbb{C}$, define
\begin{multline}
\label{definitionHfactor}
H_{\lambda/\mu}(q,p,t,b) 
:= \prod_{1\leq i < j\leq
n}\left\{\dfrac{(q^{\mu_i-\mu_{j-1}}t^{j-i})_{\mu_{j-1}-\lambda_j}
(q^{\lambda_i+\lambda_j}t^{3-j-i}b)_{\mu_{j-1}-\lambda_j}}
{(q^{\mu_i-\mu_{j-1}+1}t^{j-i-1})_{\mu_{j-1}-\lambda_j}(q^{\lambda_i
    +\lambda_j+1}t^{2-j-i}b)_{\mu_{j-1}-\lambda_j}}\right.\\
\left.\cdot 
\dfrac{(q^{\lambda_i-\mu_{j-1}+1}t^{j-i-1})_{\mu_{j-1}-\lambda_j}}
{(q^{\lambda_i-\mu_{j-1}}t^{j-i})_{\mu_{j-1}-\lambda_j}}\right\}\cdot\prod_{1\leq
i <(j-1)\leq n}
\dfrac{(q^{\mu_i+\lambda_j+1}t^{1-j-i}b)_{\mu_{j-1}-\lambda_j}}
{(q^{\mu_i+\lambda_j}t^{2-j-i}b)_{\mu_{j-1}-\lambda_j}}
\end{multline}
and 
\begin{multline}
\label{definitionSkewW}
W_{\lambda/\mu}(x; q,p,t,a,b)
:= H_{\lambda/\mu}(q,p,t,b)\cdot\dfrac{(x^{-1}, ax)_\lambda
  (qbx/t, qb/(axt))_\mu}
{(x^{-1}, ax)_\mu (qbx, qb/(ax))_\lambda}\\
\cdot\prod_{i=1}^n\left\{\dfrac{\theta(bt^{1-2i}q^{2\mu_i})}{\theta(bt^{1-2i})}
  \dfrac{(bt^{1-2i})_{\mu_i+\lambda_{i+1}}}
{(bqt^{-2i})_{\mu_i+\lambda_{i+1}}}\cdot
t^{i(\mu_i-\lambda_{i+1})}\right\}.
\end{multline}
For arbitrary $\lambda$ and $\mu$ the function
$W_{\lambda/\mu}(y, z_1, \ldots, 
z_\ell; q,p,t,a,b)$ in $\ell+1$ variables $y, z_1, \ldots, z_\ell
\in\mathbb{C}$ is defined by the following recursion formula

\begin{multline}
\label{eqWrecurrence}
W_{\lambda/\mu}(y,z_1,z_2,\ldots,z_\ell;q, p, t, a, b) \\
= \sum_{\nu\prec \lambda} W_{\lambda/\nu}(yt^{-\ell};q, p, t, at^{2\ell},
bt^\ell) \, W_{\nu/\mu}(z_1,\ldots, z_\ell;q, p, t, a, b).
\end{multline}

The definition of the \eJc\ will be needed below. Let $\lambda$ and
$\mu$ be again 
partitions of at most $n$--parts such that $\lambda/\mu$ is a skew
partition. Then the \Jc\ $\omega_{\lambda/\mu}$ are defined by
\begin{multline}
\label{eq:omega{lambda,mu}}
\omega_{\lambda/\mu}(x; r, q,p,t; a,b)
:= \dfrac{(x^{-1}, ax)_{\lambda}} {(qbx, qb/ax)_{\lambda}}
    \dfrac{(qbr^{-1}x, qb/axr)_{\mu}}{(x^{-1}, ax)_{\mu}} \\
\cdot \dfrac{(r, br^{-1}t^{1-n})_{\mu}}{(qbr^{-2}, qt^{n-1})_{\mu}}
  \prod_{i=1}^{n}\left\{ \dfrac{\theta(br^{-1}t^{2-2i} q^{2\mu_i})}
    {\theta(br^{-1}t^{2-2i})}  \left(qt^{2i-2} \right)^{\mu_i} \right\} \\
\cdot \prod_{1\leq i< j \leq n} \hspace*{-5pt} \left\{ \dfrac{
    (qt^{j-i})_{\mu_i - \mu_j} } { (qt^{j-i-1})_{\mu_i - \mu_j} }
    \dfrac{ (br^{-1}t^{3-i-j})_{\mu_i + \mu_j} } {
    (br^{-1}t^{2-i-j})_{\mu_i + \mu_j} } \right\} 
W_{\mu} (q^{\lambda}t^{\delta(n)}; q, p, t, bt^{2-2n}, br^{-1}t^{1-n})
\end{multline}
where $x,r,q,p,t, a,b\in\mathbb{C}$.

Note that $W_{\lambda/\mu}(x; q,p, t, a,b)$ vanishes unless $\lambda/\mu$
is a horizontal strip, whereas $\omega_{\lambda/\mu}(x; r;
a,b)=\omega_{\lambda/\mu}(x; r, q,p,t; a,b)$ is defined even when
$\lambda/\mu$ is not a horizontal strip. 

The operator characterization~\cite{Coskun1} of $\omega_{\lambda/\mu}$
yields a recursion formula for \Jc\ in the form 
\begin{equation}
\label{recurrence22}
  \omega_{\lambda/\tau}(y,z; r; a,b) := \sum_\mu
  \omega_{\lambda/\mu}(r^{-k}y; r; ar^{2k},
  b r^k ) \, \omega_{\mu/\tau}(z; r; a, b)
\end{equation}
where $y=(x_{1},\ldots, x_{n-k})\in\mathbb{C}^{n-k}$ and
$z=(x_{n-k+1},\ldots, x_n)\in\mathbb{C}^k$.

Using the recurrence relation~(\ref{recurrence22}) the
definition of $\omega_{\lambda/\mu}(x;r;a,b)$ can be extended 
from the single variable $x\in\mathbb{C}$ case to the multivariable
function $\omega_{\lambda/\mu}(z; r; a,b)$ with arbitrary number of variables
$z = (x_1,\ldots, x_n)\in\mathbb{C}^n$. That $\omega_{\lambda/\mu}(z;
r; a,b)$ is symmetric is also proved in~\cite{Coskun1} using a
remarkable elliptic $BC_n$ \tns\ transformation identity.

With these notation and definitions, the $M(a,b)$ and $S(b)$ matrices are now
defined. 
\begin{definition}
\label{matricesMabAndS}
Let $\lambda$ be a partition of at most $n$ parts and $q, t, a,b,\rho$
and $\sigma$ be complex parameters. Define the infinite 
matrix $M(a,b)$ by
\begin{equation}
\label{Mab}
M_{\lambda\mu}(a,b):= \dfrac{b^{|\lambda|}}{a^{|\mu|}}       
\dfrac{(a/b)_\lambda } {(qb)_\lambda }  \, K_\mu(b) \,
W_\mu(q^\lambda t^{\delta(n)}; q, p, t, at^{2-2n}, bt^{1-n}) 
\end{equation}
where 
\begin{multline}
\label{Kb}
K_{\mu}(b)=K_{\mu}(b,n):= q^{|\mu|} t^{2n(\mu)} \dfrac{
    (bt^{1-n})_{\mu}}{ (qt^{n-1})_{\mu}} \prod_{i=1}^n
  \left\{\dfrac{\theta(bt^{2-2i}q^{2\mu_i})} {\theta(b t^{2-2i})}
    \right\}  \\  \cdot  \prod_{1\leq i<j \leq 
n} \left\{\dfrac{ (qt^{j-i})_{\mu_i-\mu_j} (bt^{3-i-j})_{\mu_i+\mu_j}}
{(qt^{j-i-1})_{\mu_i-\mu_j} (b t^{2-i-j})_{\mu_i+\mu_j} }
\right\} 
\end{multline}
and the infinite diagonal matrix $S(b)$ with diagonal entries 
\begin{equation}
\label{Smatrix}
 S_{\lambda}(b) := \dfrac{(\sigma, \rho)_{\lambda}}
 { (qb/\sigma, qb/\rho)_{\lambda} }
 \left(\dfrac{qb}{\rho\sigma} \right)^{|\lambda|} 
\end{equation}
where $\abs{\lambda}=\sum_{i=1}^n \lambda_i$ and $n(\lambda) =
\sum_{i=1}^n (i-1) \lambda_i$.
\end{definition}

It will be verified that $M(a,b)$ is lower triangular in the sense defined
above, and that the $n$-dependence of definitions~(\ref{Mab})
and~(\ref{Smatrix}) is not essential. 

Let $\lambda= (\lambda_1,\lambda_2,\ldots, \lambda_n)$ be an $n$--part
partition with $\lambda_{n-m+1}\neq 0$ and $\lambda_{n-m+1} = \ldots =
\lambda_n = 0$ where $1\leq m\leq n$ and $\hat{\lambda}$ 
denote the $(n-m)$--part partition obtained by  
dropping the last $m$ zero parts from $\lambda$.  That is, 
$\hat{\lambda}=(\lambda_1,\ldots, \lambda_{n-m})$.

\begin{lemma}
\label{InvarRepLambdaM(ab)}
$M(a,b)$ is lower triangular with respect to partial inclusion
ordering, and its entries $M_{\lambda\mu}(a,b)$ are independent of any
representation of $\lambda$. That is, with the notation as above,
$M_{\lambda\mu}(a,b)$ and $M_{\hat{\lambda}\hat{\mu}}(a,b)$ are
identical.  
\end{lemma}

\begin{proof}
The fact that $M(a,b)$ is lower triangular follows from the vanishing
property~\cite{Coskun1} of $W_\mu$ function. Namely,
\begin{equation}
\label{FundamentalVanishing}
  W_{\mu}(q^\lambda t^\delta;q, p,t,a,b)=0
\end{equation}
when $\mu \nsubseteq \lambda$. 

That $M_{\lambda\mu}(a,b)$ is independent of representations of
indexing partitions follows from an analogous property~\cite{Coskun1}
of $\omega_{\lambda/\mu}$ written as
\begin{equation}
  \omega_{\lambda/\mu}(x; r; a,b) = \omega_{\hat{\lambda}/\hat{\mu}}(x; r;
  a,b).
\end{equation}
The result follows by noting that 
\begin{equation}
\label{omegaANDm}
 \omega_{\lambda\mu}(r^{-1}; r,ar^2,br) = \dfrac{(ar)_\lambda
 (qb/ar)_\mu } {(qb/a)_\lambda (ar)_\mu } \, b^{-|\lambda|+|\mu|}
 r^{|\mu|} \, M_{\lambda\mu}(br, b)          
\end{equation}
since factors of the form $u^{|\lambda|}$ and $(u)_\lambda$ are
clearly independent of $n$.
\end{proof}

A key identity in the development of the $BC_n$ \BL\ is now proved. This
identity is equivalent to the \ci~\cite{Coskun1} for
$\omega_{\lambda/\mu}$ which can be written in the form 
\begin{multline}
\label{CocycleIdentity}
  \omega_{\nu/\mu}((uv)^{-1};uv;a(uv)^2, buv) \\ = \sum_{\mu\subseteq
  \lambda \subseteq \nu} \omega_{\nu/\lambda}(v^{-1};v;a(vu)^2,bvu) \,
  \omega_{\lambda/\mu}(u^{-1};u;au^2,bu)
\end{multline}
for $u,v\in\mathbb{C}$. 

\begin{lemma}[Key Lemma]
\label{Lemma2}
With the definitions as above, 
\begin{equation}
\label{eqnLemma2}
  S^{-1}(a)\, M(c,a)\, S(a) = S^{-1}(b) \,M(c,b)\, S(b)\, M(b,a) 
\end{equation}
where $qab = c\sigma\rho$.
\end{lemma} 

\begin{proof}
The proof follows from the observation~(\ref{omegaANDm}) and the 
\ci~(\ref{CocycleIdentity}) after a simple reparametrization.
\end{proof}

Two immediate corollaries of this key result are in order.
\begin{cor}
\label{Lemma1}
For complex parameters $a,b$ and $c$,
\begin{equation}
\label{eq:Lemma1}
  M(c,a) = M(c,b)\, M(b,a) 
\end{equation}
\end{cor} 

\begin{proof}
Notice that $S(b)=I$ if $\rho=\sigma=(aq)^{1/2}$
where $I$ is the identity matrix whose entries are $I_{\lambda\mu} =
\delta_{\lambda\mu}$. 
Therefore, setting $\rho=\sigma=(aq)^{1/2}$ in the
identity~(\ref{eqnLemma2}) yields the identity~(\ref{eq:Lemma1}) to be
proved. 
\end{proof}

This result defines a cocycle relation for the matrices $M(a,b)$. The
next result shows that the matrices $M(a,b)$ are invertible just as
in the classical case~\cite{Bressoud2}. 
\begin{cor}
\label{cor:inverseM(a,b)}
For $a, b\in\mathbb{C}$,
\begin{equation}
\label{inverseM(a,b)}
  M(b, b)=I, \quad \mathrm{and} \quad M^{-1}(a,b)=M(b,a)
\end{equation}
\end{cor} 

\begin{proof}
It has been established in~\cite{Coskun1} that $\omega_{\lambda\mu}(1;
  1,a,b)=\delta_{\lambda\mu}$. Therefore the
  identity~(\ref{omegaANDm}) implies that $M(b, b)=I$. It then follows
  from Corollary~\ref{Lemma1} that $M(b,a)$ and $M(a,b)$ are
  inverses of each other.  
\end{proof}

The $M(a,b)$ matrices satisfy the following elliptic transformation
identities. 
\begin{lemma}
\label{ellipticTransfIden}
With the definitions as above,
\begin{equation}
M_{\lambda\mu}(pa, b)   
=(-1)^{|\lambda|} b^{|\lambda|-|\mu|}  a^{-|\lambda|} t^{n(\lambda)+2n(\mu)}
  q^{-|\mu|-n(\lambda')-2n(\mu')} \, M_{\lambda\mu}(a, b)   
\end{equation}
and
\begin{multline}
M_{\lambda\mu}(pa, pb) \\  
= b^{-|\lambda|} a^{|\mu|} \, (-1)^{|\lambda|+|\mu|} p^{|\lambda|-|\mu|} 
q^{-|\lambda|+|\mu|-n(\lambda')+n(\mu')} t^{n(\lambda)-n(\mu)}\,
M_{\lambda\mu}(a, b)  
\end{multline} 
\end{lemma} 

\begin{proof}
The proof follows from the observation~(\ref{omegaANDm}) and the
following identities given in~\cite{Coskun1}. 
\begin{multline}
\label{ellipticInb}
\omega_{\lambda/\mu}(x; r; a, pb) 
= (q^2 b^2 /a)^{|\lambda|-|\mu|} t^{-2n(\lambda)+2n(\mu)} 
  q^{2n(\lambda') -2n(\mu')} \omega_{\lambda/\mu}(x; r; a, b) 
\end{multline}
\begin{multline}
\label{ellipticIna}
\omega_{\lambda/\mu}(x; r; pa, b) \\ =
(q b)^{-|\lambda|+|\mu|} p^{|\lambda|-|\mu|} r^{-|\mu|}
t^{2n(\lambda)- 2n(\mu)} q^{-2n(\lambda')+2n(\mu')} 
\omega_{\lambda/\mu}(x; r; a, b) 
\end{multline}
\begin{equation}
\label{ellipticInr}
\omega_{\lambda/\mu}(x; pr; a, b) 
= (ar^{-2})^{-|\mu|} p^{2|\mu|}
  t^{2n(\mu)} q^{-2n(\mu')} \omega_{\lambda/\mu}(x; r; a, b) 
\end{equation}
where $x\in\mathbb{C}$.
\end{proof}

\subsection{Two Parameter Elliptic \BL}
The abstract matrix formulation of the two--parameter \BL\ and
consequently the notion of a \Bl\ was given in
one dimension by Agarwal, Andrews and Bressoud~\cite{AgarwalA1} and
later by Bressoud~\cite{Bressoud1}. Milne and Lilly~\cite{LillyM3},
and Andrews, Schilling and Warnaar~\cite{AndrewsS1} extended \BL\ to
root systems of type $A_n$ and $C_n$ of rank $n$ and $A_2$,
respectively~\cite{Warnaar2}. This section generalize \BL\ to elliptic
level for the non-reduced root system $BC_n$ of rank $n$, further
extending earlier results. 

The notion of a Bailey pair will be needed. Let
$\mathbb{K}$ be the field of rational functions in 
$\sigma_i,\rho_i,a_i, b_i \in\mathbb{C}$ for $i\in \mathbb{Z}_{>}$
over the field $\mathbb{C}(q,p,t)$. 
\begin{definition}
The infinite sequences $\alpha$ and $\beta$ of rational functions
$\alpha_{\lambda}, \beta_\lambda\in \mathbb{K}$  
indexed by partitions 
form a Bailey pair relative to $(b_1,a_1)$ if they satisfy   
\begin{equation}
  \beta_{\lambda} = \sum_{\mu} M_{\lambda\mu}(b_1, a_1)\,
  \alpha_{\mu}  
\end{equation}
where the sum is over partitions.
\end{definition}

The two parameter $BC_n$ \BL\ can now be proved.
\begin{theorem}[$BC_n$ \BL]
\label{twoParameterBL}
With the notation as above, suppose that $(\alpha, \beta)$ form a
Bailey pair relative to $(b_1,a_1)$. Then the pair $\beta'$ and
$\alpha'$ defined by     
\begin{equation}
  \beta'= S(a_2) S^{-1}(b_1) M(b_2,b_1) S(b_1)\, \beta
\end{equation}
and
\begin{equation}
  \alpha'= S(a_1) M(a_2,a_1)\,\alpha
\end{equation}
form a Bailey pair relative to $(b_2,a_2)$ provided that
$qa_1b_1 = a_2 \rho \sigma$.
\end{theorem} 

\begin{proof}
The proof is an immediate consequence of the Lemma~\ref{Lemma2} and
the Corollary~\ref{Lemma1}. 
\end{proof}

\begin{remark}
Note that in the iteration scheme above the parameters
$\sigma_i$, $\rho_i$ are replaced by $\sigma_{i+1}$,
$\rho_{i+1}$ and the parameters $(a_i,b_i)$ are replaced by
$(a_{i+1}, b_{i+1})$ in the $i$--th step.  
\end{remark}

\begin{remark}
Note also that the elliptic \BL\ of 
Theorem~\ref{twoParameterBL} yields a two parameter basic
(trigonometric) \BL\ when $p=0$. In that case  
$(a)_\lambda$ becomes
\begin{equation}
\label{basicQtPocSymbol}
  (a)_\lambda=(a; q, 0, t)_\lambda := \prod_{k=1}^{n} (at^{1-i};
  q,0)_{\lambda_i} = (at^{1-i}; q)_{\lambda_i}
\end{equation}
where the standard \qPs\ $(a; q)_{n}$ is as defined
in~(\ref{qPochSymbol}). The notation
$(a)_\lambda$ will be used for both the elliptic and the basic
(trigonometric) case, but the meaning will be clear from the context.  
\end{remark}

\subsection{One Parameter Basic $BC_n$ \BL}
This section and the rest of the paper uses only  
the basic (trigonometric) case $p=0$ of $BC_n$ \BL.  
Since the limiting cases of the basic $W$ functions
$W_{\lambda/\mu}(x; q,t,a,b) = W_{\lambda/\mu}(x; q, 0,t,a,b)$ will be
used in computations, some more notation is needed. Set 
\begin{equation}
\label{limit1}
  W_\mu(x; q, t, at^{2-2n}, 0) := \lim_{b
  \rightarrow 0}\, W_\mu(x; q, t, at^{2-2n}, bt^{1-n})  
\end{equation}
and
\begin{equation}
\label{limit2}
  W_\mu(x; q, t, 0 ,bt^{1-n}) := \lim_{a
  \rightarrow 0}\, (a/b)^{|\mu|}\, W_\mu(x; q, t, at^{2-2n},
  bt^{1-n})  
\end{equation}
and, finally 
\begin{equation} 
\label{limit3}
  W_\mu(q^\lambda t^{\delta(n)}; q,t, (u/v)t^{n-1})  :=
  \lim_{d\rightarrow 0} W_\mu(q^\lambda t^{\delta(n)}; q,t,
  dvt^{2-2n}, dut^{1-n}) 
\end{equation}
The existence of these limits can be seen from ($p=0$ case of) the
definition~(\ref{definitionSkewW}),  the recursion
formula~(\ref{eqWrecurrence}) and the limit rule
\begin{equation}
\label{LimitRule}
  \lim_{a\rightarrow 0}\, a^{|\mu|} (x/a)_{\mu} 
= (-1)^{|\mu|}\, x^{|\mu|} t^{-n(\mu)} q^{n(\mu')} 
\end{equation}
The limiting cases~(\ref{limit1}),~(\ref{limit2}) and~(\ref{limit3})
are essentially equivalent to known families of symmetric functions
such as Okounkov's~\cite{Okounkov1} interpolation Macdonald 
polynomials $P^*_\lambda$, $\bar{P}^*_\lambda$, which in turn
generalize Macdonald polynomials $P_\lambda$. The exact relationship
between these functions are investigated in~\cite{Coskun4}. 

Analogous to the matrix formulation of the classical
\BL~\cite{AgarwalA1}, the one parameter $BC_n$ Bailey matrix $M(b)$ is
also defined as a limiting case of the basic $M(a,b)$ matrix. 

\begin{definition}
Let $\lambda$ be a partition of at most $n$--parts and
$b\in\mathbb{C}$. Define 
\begin{equation}
\label{oneVarBM}
M_{\lambda\mu}(b):= 
L_\mu(b) \, W_\mu(q^\lambda t^{\delta(n)}; q,t, 0,
  bt^{1-n}) 
\end{equation}
where
\begin{multline}
L_\mu(b):= (-1)^{|\mu|} q^{2|\mu|+n(\mu')} t^{n(\mu)} 
\dfrac{(bt^{1-n})_\mu }{ (qt^{n-1})_\mu } \\
 \cdot  \prod_{1\leq i<j \leq 
n} \left\{\dfrac{ (qt^{j-i})_{\mu_i-\mu_j} (bt^{3-i-j})_{\mu_i+\mu_j}}
{(qt^{j-i-1})_{\mu_i-\mu_j} (b t^{2-i-j})_{\mu_i+\mu_j} }
\right\}
\end{multline}
and $n(\lambda') =\sum_{i=1}^n \binom{\lambda_i}{2}$.
\end{definition}

\begin{remark}
Note that 
\begin{equation}
M_{\lambda\mu}(b)= P_\lambda(b) \, M_{\lambda\mu}(0, b) \, Q_\mu(b)
\end{equation}
where 
\begin{equation}
P_\lambda(b) = b^{-|\lambda|} (qb)_\lambda 
\end{equation}
and
\begin{equation}
Q_\mu(b) = (-1)^{|\mu|} q^{|\mu|+n(\mu')} b^{|\mu|} t^{-n(\mu)} \,
  \prod_{i=1}^n \left\{\dfrac {(1-b
  t^{2-2i})} {(1-bt^{2-2i}q^{2\mu_i})}\right\} 
\end{equation}
Therefore it follows from the properties of the $M(a,b)$
given in Lemma~\ref{InvarRepLambdaM(ab)} that
$M(b)$ is also lower triangular and is independent of 
representations of $\lambda$. 
\end{remark}

The property that $M_{\lambda\mu}(b)$ satisfies
hyperoctahydral symmetry in the \rc\ $\lambda=k^n$ follows from the
fact that these matrix entries are well--poised hypergeometric
series. The proof of multiple $q$--series identities below uses the
symmetries that are now verified.  
\begin{lemma}
\label{M{k^nmu}(b)symmetry}
Let $\lambda=k^n$ for some non--negative integer $k$ and set $q^{z_i} =
b^{1/2} t^{1-i}$ in the definition~(\ref{oneVarBM}).
Then the matrix $M_{k^n\mu}(b)$ is invariant under the standard action
$q^{\mu+z}\leftrightarrow q^{w(\mu+z)}$ (permutations and sign
changes) of the \hgr\ of rank $n$.  
\end{lemma}

\begin{proof}
Let $x\in\mathbb{C}$. The following analogue of the \Wdegf\ for $W$
functions 
\begin{multline}
\label{conj:degform}
W_{\mu}(xt^{\delta(n)};q,t,a,b) \\
=\dfrac{(x^{-1}, axt^{n-1})_\mu}{(qbxt^{n-1}, qb/(ax))_\mu} 
\cdot \prod_{1\leq i < j\leq n} \dfrac{(t^{j-i+1})_{\mu_i
-\mu_j}(qbt^{n-i-j+1})_{\mu_i+\mu_j}}
{(t^{j-i})_{\mu_i -\mu_j}(qbt^{n-i-j})_{\mu_i+\mu_j}}
\end{multline}
follows from the basic (i.e., $p=0$) version of the
$W$--\Js~\cite{Coskun1}    
\begin{multline}
\label{eq:W_Jackson}
W_{\lambda}(s^{-1}x;q,p,t,at^{-2n}s^2,bt^{-n}s)\\
=\dfrac{(s , ast^{-n-1})_{\lambda}}
{(qbt^{-1}, qbt^n/a)_{\lambda}}
\cdot \prod_{1\leq i < j\leq n}
\left\{\dfrac{(t^{j-i+1})_{\lambda_i-\lambda_j}
(qbst^{-i-j+1})_{\lambda_i+\lambda_j}}
{(t^{j-i})_{\lambda_i-\lambda_j} (qbst^{-i-j})_{\lambda_i +
    \lambda_j}} \right\}  \\
\cdot\sum_{\mu \subseteq \lambda}
  \dfrac{(bt^{-n}, qbt^n/(as))_\mu}{(qt^{n-1}, 
ast^{-n-1})_\mu} \cdot \prod_{i=1}^n \left\{
  \dfrac{\theta(bt^{1-2i}q^{2\mu_i})}{\theta(bt^{1-2i})}(qt^{2i-2})^{\mu_i}
\right\} \\ \cdot \prod_{1\leq i < j \leq n}
\left\{ \dfrac{(t^{j-i})_{\mu_i -\mu_j} (qt^{j-i})_{\mu_i
      -\mu_j}}{(qt^{j-i-1})_{\mu_i -\mu_j}(t^{j-i+1})_{\mu_i
      -\mu_j}} \dfrac{(bqt^{-i-j})_{\mu_i+\mu_j}
    (bt^{-i-j+2})_{\mu_i+\mu_j}} {(bt^{-i-j+1})_{\mu_i+\mu_j}
(qbt^{-i-j+1})_{\mu_i+\mu_j}} \right\} \\
\cdot W_\mu(q^\lambda t^{\delta(n)};q,p,t,bst^{1-2n},bt^{-n})\cdot
W_\mu(x;q,p,t,at^{-2n},bt^{-n})
\end{multline}
and the fact~\cite{Coskun1} that 
\begin{equation}
W_{\mu}(t^{\delta(n)};q,t,a,b) = \delta_{0\mu}.
\end{equation}

Set $\lambda=k^n$ and $q^{z_i} = b^{1/2} t^{1-i}$
in~(\ref{oneVarBM}). Using the identity~(\ref{conj:degform}), the
definition of the \qPs~(\ref{qPochSymbol}), and the identities
\begin{equation}
\label{qPsFlipIdentity}
  (v)_n = \dfrac{(-v)^n q^{\binom{n}{2}}} { (qv^{-1})_{-n}}
\end{equation}
and
\begin{equation}
\label{doubleProductOne}
\prod_{1\leq i<j \leq n} \dfrac{a_i}{a_j} 
= \prod_{ i=1 }^{n} a_i^{1+n-2i} 
\end{equation}
the entries of the matrix $M(b)$ can be written in the form
\begin{multline}
\label{M(b)matrixIsInvariant}
M_{k^n\mu}(b)= 
 \prod_{i=1}^n \dfrac{(q^{-k}b^{-1/2}q^{-z_i})_{\infty}  
  (q^{-k}b^{-1/2}q^{z_i})_{\infty}  }
  {(qb^{-1/2}t^{n-1}q^{-z_i})_{\infty}
 (qb^{-1/2}t^{n-1}q^{z_i})_{\infty} }  \\ \cdot 
\prod_{1\leq i<j \leq n} \left\{\dfrac{ (q^{1-z_j+z_i})_{\infty}
  (q^{1+z_j-z_i})_{\infty} } 
{(qt^{-1}q^{z_i-z_j} )_{\infty} (qt^{-1} q^{-z_i+z_j})_{\infty} }
   \dfrac{(q^{1-z_i-z_j})_{\infty} (q^{1+z_i+z_j} )_{\infty}}
 {(t^{-1} q^{1-z_i-z_j})_{\infty} (t^{-1} q^{1+z_i+z_j} )_{\infty}}
 \right\} \\ \cdot 
\prod_{i=1}^n  \dfrac{(qb^{-1/2}t^{n-1}q^{-\mu_i-z_i} )_{\infty}  
 (qb^{-1/2}t^{n-1}q^{\mu_i+z_i} )_{\infty}  }
 {(q^{-k}b^{-1/2}q^{-\mu_i-z_i} )_{\infty}  
  (q^{-k}b^{-1/2}q^{\mu_i+z_i} )_{\infty} } \\ \cdot
\prod_{1\leq i<j \leq n} \left\{\dfrac{ (q^{1+\mu_i+z_i-\mu_j-z_j})_{\infty}
  (q^{1-\mu_i-z_i +\mu_j+z_j})_{\infty}} 
{(t^{-1}q^{1+\mu_i+z_i-\mu_j-z_j})_{\infty}
 (t^{-1}q^{1-\mu_i-z_i +\mu_j+z_j })_{\infty}} \right. \\ \left.
  \cdot \dfrac{(q^{1-\mu_i-z_i-\mu_j-z_j})_{\infty}
 (q^{1+\mu_i+z_i+\mu_j+z_j})_{\infty}} 
 {(t^{-1}q^{1-\mu_i-z_i-\mu_j-z_j})_{\infty}
  (t^{-1}q^{1+\mu_i+z_i+\mu_j+z_j})_{\infty}} \right\}
\end{multline}
It is now clear that $M_{k^n\mu}(b)$ is invariant under the maps
$q^{\mu_i+z_i} \leftrightarrow q^{\pm (\mu_j+z_j)}$. 
\end{proof}

It could be shown, using the definition and the recurrence relation
for $W_\lambda$, that $M_{\lambda\mu}(b)$ satisfies \hgr\
symmetry for a general partition $\lambda$. However, the \rc\ given in
Lemma~\ref{M{k^nmu}(b)symmetry} will be sufficient for the proof of
the multiple $q$--series identities obtained below. 

The diagonal $S(b)$ matrix has similar properties. 
\begin{lemma}
\label{S{lambda}(b)symmetry}
The entries of the basic $S(b)$ matrix~(\ref{Smatrix}) is independent of
representations of $\lambda$ and satisfies the \hgr\ symmetry. 
\end{lemma}

\begin{proof}
That $S_\lambda(b)$ is independent of different representations of
$\lambda$ is obvious. Setting $q^{z_i} = b^{1/2} t^{1-i}$
in~(\ref{Smatrix}) and using the identities given above, one obtains 
\begin{multline}
 S_\lambda(b)
= \prod_{i=1}^n \dfrac{(\sigma b^{-1/2} q^{z_i})_{\infty}
 (\sigma b^{-1/2} q^{-z_i})_{\infty} }
 { (q\rho^{-1} b^{1/2} q^{z_i} )_{\infty} (q \rho^{-1} b^{1/2}
   q^{-z_i} )_{\infty}} \\ \cdot 
\dfrac{ (q\rho^{-1} b^{1/2} q^{z_i+\lambda_i} )_{\infty} (q \rho^{-1}
   b^{1/2} q^{-(z_i+\lambda_i)} )_{\infty}} {(\sigma b^{-1/2}
   q^{z_i+\lambda_i})_{\infty} 
 (\sigma b^{-1/2} q^{-(z_i+\lambda_i)})_{\infty} }
\end{multline}
proving the invariance under the maps
$q^{\lambda+z} \leftrightarrow q^{w(\lambda+z)}$ or $q^{\lambda_i+z_i} 
\leftrightarrow q^{\pm(\lambda_i+z_i)}$ for all $i\in [n]$. 
\end{proof}

The fact that $M(b)$ is invertible is a simple consequence of
the Corollary~\ref{Lemma1}. 
\begin{lemma}
\label{inverseOfM}
With the notation as above, the inverse $M^{-1}(b)$ of the infinite
triangular matrix $M(b)$ is given by 
\begin{multline}
  M^{-1}_{\lambda\mu}(b)
  = \dfrac{q^{-|\lambda| + |\mu|} t^{2n(\mu)} }{(qb, qt^{n-1})_\mu}
  \prod_{i=1}^n \left\{\dfrac{(1-bt^{2-2i}q^{2\lambda_i})} {(1-b
  t^{2-2i})} \right\} \\ \cdot \prod_{1\leq
  i<j \leq  n} \left\{\dfrac{ (qt^{j-i})_{\mu_i-\mu_j} }
{(qt^{j-i-1})_{\mu_i-\mu_j} } \right\} \cdot W_\mu(q^\lambda
  t^{\delta(n)}; q,t, bt^{2-2n}, 0) 
\end{multline}
\end{lemma}

\begin{proof}
Setting $c=a$ in Lemma~(\ref{Lemma1}) gives the identity $I=M(a,b)\,M(b,a)$
by the virtue of~(\ref{inverseM(a,b)}). This beautiful identity can be
written explicitly in the form 
\begin{multline}
\label{orthogonalityForW}
\delta_{\lambda\tau} = \sum_\mu \dfrac{b^{|\lambda|}}{a^{|\mu|}}       
\dfrac{(a/b)_\lambda } {(qb)_\lambda } \, K_\mu(b) \, W_\mu(q^\lambda
t^{\delta(n)}; q,t, at^{2-2n}, bt^{1-n}) \\
\cdot \dfrac{a^{|\mu|}}{b^{|\tau|}}       
\dfrac{(b/a)_\mu } {(qa)_\mu } \, K_\tau(a) \, W_\tau(q^\mu t^{\delta(n)};
  q,t, bt^{2-2n}, at^{1-n})
\end{multline}
Noting that 
\begin{equation}
K_{\mu}(b)=  (-1)^{|\mu|} q^{-|\mu|-n(\mu')} 
  t^{n(\mu)} \prod_{i=1}^n
  \left\{\dfrac{(1-bt^{2-2i}q^{2\mu_i})} {(1-b t^{2-2i})}
    \right\} \, L_\mu(b)
\end{equation}
and taking the limit $a\rightarrow 0$ of~(\ref{orthogonalityForW}) 
by using~(\ref{LimitRule}), one gets the inverse of $M(b)$ given
above. 
\end{proof}

Next, the closed form of the entries of the matrix $N(b)
:= M(b)\, S(b)\, M^{-1}(b)$ will be computed. In the classical \od\
case, $N(b)$ is computed by means of 
the \BT\ and the Saalchutz's \tts\ summation formula. Instead of
generalizing this idea, the
two parameter basic $BC_n$ \BL\ is used to prove the next result. 
\begin{lemma}
Let $q, t, b, \sigma, \rho \in\mathbb{C}$. Then the entries of the
matrix $N(b)$ is given in the closed form 
\begin{multline}
\label{N(b)matrix}
  N_{\lambda\mu}(b)
= q^{|\mu|} t^{2n(\mu)} \dfrac{ (qb)_\lambda
  (qb/\rho\sigma)_\lambda}{(qb/\sigma)_{\lambda} (qb/\rho)_{\lambda} }  
\dfrac{ (\sigma)_{\mu} (\rho)_{\mu} }{(qb)_\mu (qt^{n-1})_{\mu} } \\
\prod_{1\leq i<j \leq n} \left\{\dfrac{ (qt^{j-i})_{\mu_i-\mu_j} }
{(qt^{j-i-1})_{\mu_i-\mu_j} }
\right\} \cdot W_\mu(q^\lambda t^{\delta(n)}; q,t, qb t^{n-1}/\rho\sigma) 
\end{multline}
\end{lemma}

\begin{proof}
Note that one can write 
\begin{equation}
M_{\lambda\mu}(b) = P_\lambda(b) \, M_{\lambda\mu}(0, b) \, Q_\mu(b)
\end{equation}
where 
\begin{equation}
P_\lambda(b) = b^{-|\lambda|} (qb)_\lambda 
\end{equation}
and
\begin{equation}
Q_\mu(b) = (-1)^{|\mu|} q^{|\mu|+n(\mu')} b^{|\mu|} t^{-n(\mu)} \,
  \prod_{i=1}^n \left\{\dfrac {(1-b
  t^{2-2i})} {(1-bt^{2-2i}q^{2\mu_i})}\right\} 
\end{equation}
Set $a=du$ and $c=dv$ for some $d,u,v\in\mathbb{C}$
in~(\ref{eqnLemma2}). Due to the relation $qub=v\rho\sigma$ 
one gets $v/u=qb/\rho\sigma$. Rewrite Lemma~(\ref{Lemma2}) in the form  
\begin{equation}
  M(dv,b) \, S(b) \, M(b,du)  
= S(b)\, S^{-1}(du) \, M(dv,du) \, S(du)  
\end{equation}
Now in the limit as $d\rightarrow 0$, the \lhs\ becomes
\begin{equation}
\lim_{d\rightarrow 0} M(dv,b) \, S(b) \, M(b,du)  
= M(0,b) \, S(b) \, M(b,0)  
\end{equation}
which is essentially equal to $N(b)$ up to diagonal factors. That is,
\begin{equation}
\label{N(b)}
  N(b) = \lim_{d\rightarrow 0} \left\{ P \,S(b)\, S^{-1}(du) \,
  M(dv,du) \, S(du) \, P^{-1}\right\} 
\end{equation}
which gives the closed form to be computed. 
\end{proof}

The one parameter basic $BC_n$ \BL\ is proved next. The
notion of a \Bp\ translates to the one parameter case in an obvious
way. Namely, the infinite sequences $\alpha$ and $\beta$ of rational
functions $\alpha_{\lambda}, \beta_\lambda\in \mathbb{K}$  
indexed by partitions form a Bailey pair (relative to $b=b_1$) if they
satisfy $\beta_{\lambda} = \sum_{\mu} M_{\lambda\mu}(b)\,
\alpha_{\mu}$ where the sum is over partitions.

\begin{theorem}
\label{OneParaBaileyLem}
Suppose that the infinite sequences $\alpha$ and $\beta$ form a Bailey
pair relative to $b$. Then $\alpha'$ and $\beta'$ also form
a Bailey pair relative to $b$ where  
\begin{equation}
  \alpha'_{\lambda} = S_{\lambda}(b)\, \alpha_{\lambda} 
\end{equation}
and
\begin{equation}
  \beta'_{\lambda} = \sum_{\mu} N_{\lambda\mu}(b)\, \beta_{\mu} 
\end{equation}
where the sum is over partitions.
\end{theorem}

\begin{proof}
The proof is a consequence of the definition of $N(b)$ and that of a
\Bp\ relative to $b$.
\end{proof}

\subsection{Generalized \Wt} 
\label{subsection:Iteration}
The power of \BL\ comes from its potential for 
iteration. The lemma can be 
applied to a given \Bp\ $(\alpha, \beta)$ repeatedly producing an
infinite sequence of \Bp s
$(\alpha,\beta)\rightarrow(\alpha',\beta')\rightarrow  
(\alpha'',\beta'')\rightarrow\cdots$, what is called a \Bc. In fact,
a stronger result says that it is possible to walk along the \Bc\
in every direction as depicted in the following figure.  
\begin{figure}[h]
\label{BaileyChain}
\begin{displaymath}
\begin{array}{ccccccccccccc}
\cdots & \leftrightarrow & \alpha^{(-2)} & \leftrightarrow &
\alpha^{(-1)} & \leftrightarrow & \alpha & \leftrightarrow &
\alpha' & \leftrightarrow & \alpha'' & \leftrightarrow & \cdots \\
 & & \updownarrow & & \updownarrow & & \updownarrow & &
\updownarrow & & \updownarrow &  \\
\cdots & \leftrightarrow & \beta^{(-2)} & \leftrightarrow &
\beta^{(-1)} & \leftrightarrow & \beta & \leftrightarrow &
\beta' & \leftrightarrow & \beta'' & \leftrightarrow & \cdots
\end{array}
\end{displaymath}
\caption{\Bc}
\end{figure}

\begin{lemma}[Bailey Walk]
The entire \Bc\ is uniquely determined by a single node $\alpha^{(i)}$
or $\beta^{(i)}$ for any $i\in\mathbb{Z}$ on the chain.
\end{lemma}

\begin{proof}
The lower triangular matrices $M(b)$, $N(b)$ and the diagonal
matrix $S(b)$, having no zero entries on their diagonal, are all
invertible. One move forward and backward in the first line in
Figure 1 by $S(b)$
and $S^{-1}(b)$, in the second line by $N(b)$ and $N^{-1}(b)$ and move
up and down by $M(b)$ and $M^{-1}(b)$. $M^{-1}(b)$ is obtained 
before in Lemma~(\ref{inverseOfM}), and the fact that $N(b)$ is
invertible follows from its definition. 
\end{proof}

This powerful iteration mechanism allows one to prove numerous
multiple \bhs\ and multiple $q$-series identities. In this section
a terminating \sfs\ summation formula and a generalized
\Wt\ will be given. The limiting cases of these results are used to
prove \Epnt, the \RRis\ and the extreme cases of the \AGis. The
details for the \RRis\ will be given in this paper and other results
will appear in future publications~\cite{Coskun3},~\cite{Coskun4}.

The \Bp\ $(\alpha, \beta)$ corresponding to the simplest non--trivial
sequence $\beta$ defined by $\beta_{\lambda} = \delta_{\lambda 0}$ 
is called the unit \Bp. 
The corresponding $\alpha$ sequence can easily be computed using the
the inverse matrix $M^{-1}(b)$. One gets
\begin{equation}
\label{UBPalpha}
  \alpha_{\lambda} = \sum_{\mu} M^{-1}_{\lambda\mu}(b) \,\beta_{\mu} =
  q^{-|\lambda| } \prod_{i=1}^n
  \left\{\dfrac{(1-bt^{2-2i}q^{2\lambda_i})} {(1-b t^{2-2i})} \right\}
\end{equation}
Iterating the \BL\ of Theorem~\ref{OneParaBaileyLem} once, that is
computing the sequences $\alpha'$ 
and $\beta'$ and writing out the relation $\beta'_{\lambda} =
\sum_{\mu} M_{\lambda\mu}(b)\, \alpha'_{\mu} $ explicitly gives a
\hd\ analogue of the terminating $_6\varphi_5$ summation formula
\begin{multline}
\label{sfs}
\dfrac{ (qb, qb/\rho_1\sigma_1)_\lambda}{(qb/\sigma_1,
  qb/\rho_1)_{\lambda} } =  \sum_{\mu\subseteq \lambda} (-1)^{|\mu|}
  q^{n(\mu')} t^{n(\mu)} \left( \dfrac{q^2b}{\sigma_1\rho_1}
  \right)^{|\mu|} \dfrac{(bt^{1-n}, \sigma_1,
  \rho_1)_{\mu}}{(qt^{n-1}, qb/\sigma_1, qb/\rho_1)_{\mu}} 
\\ \cdot  \prod_{i=1}^n
  \left\{\dfrac{(1-bt^{2-2i}q^{2\mu_i})} {(1-b t^{2-2i})} \right\}
  \prod_{1\leq i<j \leq n} \left\{\dfrac{ (qt^{j-i})_{\mu_i-\mu_j}
  (bt^{3-i-j})_{\mu_i+\mu_j}} {(qt^{j-i-1})_{\mu_i-\mu_j} (b
  t^{2-i-j})_{\mu_i+\mu_j} } \right\} \\ \cdot  W_\mu(q^\lambda
  t^{\delta(n)}; q,t, 0, bt^{1-n})  
\end{multline}
A second iteration of \BL\ and
consequently writing up the relation $\beta^{''}_{\lambda} =
\sum_{\mu} M_{\lambda\mu}(b)\, \alpha^{''}_{\mu}$ gives
\begin{multline}
\label{WatsonTransform}
 \dfrac{ (qb, qb/\rho_2\sigma_2)_\lambda}{(qb/\sigma_2,
  qb/\rho_2)_{\lambda} } \sum_{\mu\subseteq \lambda} q^{|\mu|}
  t^{2n(\mu)} \dfrac{ (\sigma_2, \rho_2, qb/\rho_1\sigma_1)_\mu} {(qt^{n-1},
  qb/\sigma_1, qb/\rho_1)_{\mu} } \\
  \cdot \prod_{1\leq i<j \leq n} \left\{\dfrac{
  (qt^{j-i})_{\mu_i-\mu_j} } {(qt^{j-i-1})_{\mu_i-\mu_j} } \right\} \,
  W_\mu(q^\lambda t^{\delta(n)}; q,t, qbt^{n-1}/\rho_2\sigma_2 )
   \\ = \sum_{\mu\subseteq \lambda} \left(
  \dfrac{q^3b^2}{\sigma_1\rho_1\sigma_2\rho_2} \right)^{|\mu|}
  (-1)^{|\mu|} q^{n(\mu')} t^{n(\mu)} \prod_{i=1}^n
  \left\{\dfrac{(1-bt^{2-2i}q^{2\mu_i})} {(1-b t^{2-2i})} \right\}
\\ \cdot \prod_{1\leq i<j
  \leq n} \left\{\dfrac{ (qt^{j-i})_{\mu_i-\mu_j}
  (bt^{3-i-j})_{\mu_i+\mu_j}} {(qt^{j-i-1})_{\mu_i-\mu_j} (b
  t^{2-i-j})_{\mu_i+\mu_j} } \right\}\, W_\mu(q^\lambda t^{\delta(n)};
  q,t, 0, bt^{1-n}) \\ \cdot \dfrac{(bt^{1-n}, \sigma_2, \rho_2, \sigma_1,
  \rho_1)_{\mu}}{(qt^{n-1}, qb/\sigma_1, qb/\rho_1, qb/\sigma_2,
  qb/\rho_2)_{\mu}} 
\end{multline}
This is a \hd\ analogue of the \Wt~(\ref{eq:onedimWt}). 

Iterating the \BL\ of Theorem~\ref{OneParaBaileyLem} $N$
  times yields an extension of the \Wt~(\ref{WatsonTransform}) which
  is called the generalized \Wt. The following notation will be used. 
\begin{equation}
  \sum_{{\nu_{n}}\subseteq{\nu_{n-1}}\subseteq \cdots \subseteq
  {\nu_1} \subseteq\nu_{0} } :=
\sum_{{\nu_1}\subseteq\nu_{0}} \sum_{{\nu_2}\subseteq{\nu_1}}  \cdots
  \sum_{{\nu_{n}}\subseteq{\nu_{n-1}}} 
\end{equation}

\begin{lemma}
\label{lem:generalWatson}
With the notation as above,
\begin{multline}
\label{generalWatson}
\sum_{\mu^{N-1}
    \subseteq \cdots \subseteq  \mu^0} 
    \!\! \prod_{k=1}^{N-1} \left\{ q^{|\mu^k|} t^{2n(\mu^k)} \dfrac{
    (qb/\rho_{N-k} 
    \sigma_{N-k}, \sigma_{N-k+1}, \rho_{N-k+1})_{\mu^{k}} }
    {(qb/\sigma_{N-k}, qb/\rho_{N-k}, qt^{n-1})_{\mu^k} }  
 \right. \\ \left. \cdot \prod_{1\leq i<j \leq n} \left\{\dfrac{
    (qt^{j-i})_{\mu^k_i-\mu^k_j} } {(qt^{j-i-1})_{\mu^k_i-\mu^k_j} }
    \right\} \, W_{\mu^k}(q^{\mu^{k-1}} t^{\delta(n)}; q,t, 
    qbt^{n-1}/\rho_{N-k+1}\sigma_{N-k+1} ) \right\} \\ =
    \dfrac{(qb/\sigma_N, qb/\rho_N)_{\mu^0} } { (qb,
    qb/\rho_N\sigma_N)_{\mu^0} } \cdot \!\!\!  
    \sum_{\mu\subseteq \mu^0} 
    (-1)^{|\mu|} q^{|\mu|+n(\mu')} t^{n(\mu)} \dfrac{(bt^{1-n})_\mu }{
    (qt^{n-1})_\mu } \prod_{i=1}^n
    \left\{\dfrac{(1-bt^{2-2i}q^{2\mu_i})} {(1-b t^{2-2i})} \right\}
    \\ \cdot  \prod_{1\leq i<j \leq n} \left\{\dfrac{
    (qt^{j-i})_{\mu_i-\mu_j} (bt^{3-i-j})_{\mu_i+\mu_j}}
    {(qt^{j-i-1})_{\mu_i-\mu_j} (b t^{2-i-j})_{\mu_i+\mu_j} }
    \right\}\, W_\mu(q^{\mu^0} t^{\delta(n)}; q,t, 0, bt^{1-n}) \\
    \cdot \prod_{k=1}^N \left\{ \dfrac{(\sigma_{N-k+1},
    \rho_{N-k+1})_{\mu} } {(qb/\sigma_{N-k+1}, qb/\rho_{N-k+1})_{\mu}
    } \left( \dfrac{qb}{\sigma_{N-k+1}\rho_{N-k+1}} \right)^{|\mu|}
    \right\} 
\end{multline}
\end{lemma}

\begin{proof}
Iterate the \BL\ of Theorem~\ref{OneParaBaileyLem} $N$ times
starting with the unit \Bp\ corresponding to $\beta_{\lambda} =
\delta_{\lambda 0}$.   
\end{proof}

Note that the $N=1$ case of the generalized \Wt~(\ref{generalWatson})
reduces to the terminating \sfs\ summation formula~(\ref{sfs}), and the
$N=2$ case of it reduces to the \Wt~(\ref{WatsonTransform}). 

\begin{remark}
\label{HyperoctahedralSymmetry}
It has been already seen above in Lemma~\ref{M{k^nmu}(b)symmetry} and
Lemma~\ref{S{lambda}(b)symmetry} that the matrices $M(b)$ and $S(b)$
involved in the iteration 
process are invariant under the action of the \hgr, and are
independent of different representations of $\lambda$. Therefore, the
\rhs\ of the identity~(\ref{generalWatson}) obtained in the iteration
process will always share the same properties when 
$\alpha_\lambda$ has the same properties. But, if one sets $q^{z_i} =
b^{1/2} t^{1-i}$ in~(\ref{UBPalpha}), $\alpha_{\lambda}$ can be written 
in the form 
\begin{equation}
  \alpha_{\lambda} = 
 \prod_{i=1}^n \left\{ q^{-2 (z_i+\lambda_i)^2 +2z_i^2 } 
   \dfrac{(q^{1+2z_i} )_{\infty}}{(q^{1+2(z_i+\lambda_i)})_\infty}
   \dfrac{(q^{1-2z_i} )_{\infty}}{(q^{1-2(z_i+\lambda_i)})_\infty}
 \right\} 
\end{equation}
which clearly satisfies these properties. To make the symmetries more
transparent, one 
can rewrite the series in the \rhs\ of the generalized
\Wt~(\ref{generalWatson}) similar to~(\ref{M(b)matrixIsInvariant}) by
setting $\mu^0=k^n$ and $q^{z_i} = 
b^{1/2} t^{1-i}$ and using the same definitions and identities, namely
(\ref{qPochSymbol}), (\ref{qPsFlipIdentity}),
and~(\ref{doubleProductOne}) as follows. 
\begin{multline}
\label{genRSiHyperSym}
 \prod_{i=1}^n q^{-(N-1)z_i^2 } \prod_{i=1}^{n}
   \dfrac{(q^{1+2z_i})_\infty (q^{1-2z_i})_\infty} 
   {(qt^{n-i})_\infty (q^{1-2z_i}t^{n-i})_\infty } \\  
 \cdot \prod_{1\leq i<j \leq n} \left\{ \dfrac{ (q^{1+z_i-z_j})_\infty
   (q^{1-z_i+z_j})_\infty } 
 {(t^{-1}q^{1+z_i-z_j})_\infty (t^{-1}q^{1-z_i+z_j})_\infty } 
   \dfrac{ (q^{1-z_i-z_j})_\infty
   (q^{1+z_i+z_j})_\infty } 
 {(t^{-1}q^{1-z_i-z_j})_\infty (t^{-1}q^{1+z_i+z_j})_\infty } \right\} \\
 \cdot \sum_{\mu,\ell(\mu)\leq n} \prod_{i=1}^n
   q^{(N-1)(\mu_i+z_i)^2 }  \prod_{i=1}^{n}
   \dfrac{(q^{1-z_i+(\mu_i+z_i)}t^{n-i})_\infty
   (q^{1-z_i-(\mu_i+z_i)}t^{n-i})_\infty 
   } {(q^{1+2z_i+2\mu_i})_\infty (q^{1-2z_i-2\mu_i})_\infty}  \\  
 \cdot \prod_{1\leq i<j \leq n} \left\{\dfrac{
   (q^{1+z_i-z_j+\mu_i-\mu_j})_\infty 
   (q^{1-z_i+z_j-\mu_i +\mu_j})_\infty } 
 {(t^{-1}q^{1+z_i-z_j+\mu_i-\mu_j})_\infty (t^{-1}q^{1-z_i+z_j-\mu_i
   +\mu_j})_\infty }  \right. \\ \left.
 \cdot \dfrac{ (q^{1-z_i-z_j-\mu_i-\mu_j})_\infty 
   (q^{1+z_i+z_j+\mu_i +\mu_j})_\infty } 
 {(t^{-1}q^{1-z_i-z_j-\mu_i-\mu_j})_\infty (t^{-1}q^{1+z_i+z_j+\mu_i
   +\mu_j})_\infty }  \right\}
\end{multline}
\end{remark}

\section{Multiple $q$--Series Identities Associated to Root Systems}
\label{section4}

The generalized \Wt~(\ref{generalWatson}) produces, in the limit,
several remarkable multiple $q$--series identities. The initial $N=1$
instance of the iteration corresponds to the \Epnt, and the $N=2$
instance yields the \RRis. Furthermore, the general $N$ case
of~(\ref{generalWatson}) is enough to prove the
\AGis\ in the extreme cases~\cite{Coskun1}. However the full \AGis\
requires the two parameter $BC_n$ \BL\ of
Theorem~\ref{twoParameterBL}. As noted above, only the details
for the \RRis\ will be given in this paper.

\subsection{The Rogers--Selberg Identity} 
First, a $BC_n$ analogue of the \RSi~(\ref{eq:onedimRSi}) will be
proved. 
\begin{lemma}
\label{lem:unspecializedRSi}
Let $q, t, b \in \mathbb{C}$ and $\abs{q}<1$. A $BC_n$
generalization of the \RSi\ is given by 
\begin{multline}
\label{unspecializedRSi}
\sum_{\ell(\lambda)\leq n}  (-1)^{|\lambda|}\, 
b^{2|\lambda|} t^{(1-n)|\lambda|-3n(\lambda)}
q^{2|\lambda|+5n(\lambda')} \dfrac{   
    (bt^{1-n})_{\lambda}}{ (qt^{n-1})_{\lambda}}       
  \prod_{i=1}^{n} \left\{ \dfrac{(1-bt^{2-2i} q^{2\lambda_i})}
    {(1-bt^{2-2i})}  \right\} \\ \cdot
  \prod_{1\leq i< j \leq n} \left\{ \dfrac{ (qt^{j-i})_{\lambda_i - \lambda_j}
    } { (qt^{j-i-1})_{\lambda_i - \lambda_j} } \dfrac{
    (bt^{3-i-j})_{\lambda_i + \lambda_j} } { (bt^{2-i-j})_{\lambda_i + 
    \lambda_j} } \dfrac{(t^{j-i+1})_{\lambda_i
-\lambda_j}(qbt^{2-i-j})_{\lambda_i+\lambda_j}}
{(t^{j-i})_{\lambda_i -\lambda_j}(qbt^{1-i-j})_{\lambda_i+\lambda_j}}
\right\}   \\ 
= (qb)_{\infty^n} \!\! \sum_{\ell(\lambda)\leq n} \dfrac{  
b^{|\lambda|} q^{|\lambda|+2n(\lambda')} t^{(1-n)|\lambda|} }{
(qt^{n-1})_{\lambda}}  \!\!\!\!
  \prod_{1\leq i< j \leq n} \left\{ \dfrac{ (qt^{j-i})_{\lambda_i -
    \lambda_j} } { (qt^{j-i-1})_{\lambda_i - \lambda_j} }
\dfrac{(t^{j-i+1})_{\lambda_i -\lambda_j} } {(t^{j-i})_{\lambda_i
    -\lambda_j} }   \right\} 
\end{multline}
where $(u)_{\infty^n}$ denotes the product $\prod_{i=1}^n
(ut^{1-i})_\infty $. 
\end{lemma}

\begin{proof}
Set $\lambda=k^n$, for some $k\in\mathbb{Z}_{\geq}$, in the
\Wt~(\ref{WatsonTransform}) 
and send the parameters $\sigma_1, \rho_1, \sigma_2$ and $\rho_2$ to
$\infty$ using the limit rule~(\ref{LimitRule}) to get 
\begin{multline}
\label{limitingWatson}
\dfrac{1 } {(qb)_{k^n} } 
\sum_{\lambda\subseteq k^n} 
  t^{(1-n)|\lambda|-2n(\lambda)} \,b^{2|\lambda|}
  \,q^{(2+k)|\lambda|+4n(\lambda')}   \\ \cdot
  \prod_{i=1}^{n}\left\{ \dfrac{(1-bt^{2-2i} q^{2\lambda_i})}
    {(1-bt^{2-2i})}  \right\}  \dfrac{ 
    (bt^{1-n})_{\lambda}}{ (qt^{n-1})_{\lambda}} \,
  \dfrac{(q^{-k})_\lambda} {(q^{1+k}b)_\lambda } 
    \\ \cdot  \prod_{1\leq i< j \leq n} \left\{ \dfrac{
        (qt^{j-i})_{\lambda_i - \lambda_j} 
    } { (qt^{j-i-1})_{\lambda_i - \lambda_j} } \dfrac{
    (bt^{3-i-j})_{\lambda_i + \lambda_j} } { (bt^{2-i-j})_{\lambda_i + 
    \lambda_j} } \dfrac{(t^{j-i+1})_{\lambda_i
-\lambda_j}(qbt^{2-i-j})_{\lambda_i+\lambda_j}}
{(t^{j-i})_{\lambda_i -\lambda_j}(qbt^{1-i-j})_{\lambda_i+\lambda_j}}
\right\} \\
=  \sum_{\lambda\subseteq k^n}
      q^{(k+1)|\lambda|+ n(\lambda')} \, (-1)^{|\lambda|}\,
      b^{|\lambda|}\, t^{(1-n)|\lambda|+n(\lambda)} 
\dfrac{ (q^{-k})_\lambda }{ (qt^{n-1})_{\lambda}} \,
 \\ \cdot \prod_{1\leq i< j \leq n} \left\{ \dfrac{ (qt^{j-i})_{\lambda_i -
    \lambda_j} } { (qt^{j-i-1})_{\lambda_i - \lambda_j} }
\dfrac{(t^{j-i+1})_{\lambda_i -\lambda_j} } {(t^{j-i})_{\lambda_i
    -\lambda_j} }   \right\} 
\end{multline}

Before passing the limit $k\rightarrow \infty$, it would be useful to
rewrite this limiting case of the \Wt~(\ref{limitingWatson}) in order
to verify convergence. Note that factors in the
series~(\ref{limitingWatson}) may be flipped using 
\begin{equation}
  (v;q,t)_\mu = \dfrac{1} { (qt^{n-1}/v; q, t)_{-\mu^{(r)}} }
  (-1)^{|\mu|} v^{|\mu|} q^{n(\mu')} \, t^{-n(\mu)} 
\end{equation}
where $-\mu^{(r)}$ is defined to be 
\begin{equation}
\label{reversedPart}
  -\mu^{(r)} := (-\mu_n,\ldots,-\mu_{n+1-i},\ldots, -\mu_1)
\end{equation}
Set $q^{z_i}=b^{1/2} t^{1-i}$ in~(\ref{limitingWatson}), flip
appropriate factors and use the
definition of the \qPs~(\ref{qPochSymbol}) to write the \wps\ (i.e.,
\lhs) and the \bs\ (i.e., \rhs) of~(\ref{limitingWatson}) in the form  
\begin{multline}
\label{wpsWt2ndversion}
\prod_{i=1}^n \dfrac{ (q^{2z_i}t^{i-n})_{\infty} } { (qt^{n-i})_{\infty} }
\dfrac{(q^{1+2z_i} )_{\infty}} { (q^{2z_i})_{\infty} } \dfrac{1}
  {(q^{1+k} t^{i-1} )_\infty (q^{1+k+2z_i} t^{i-1} )_\infty } \\ \cdot
  \prod_{1\leq i< j \leq n} \left\{ \dfrac{ (q^{1+z_i-z_j})_{\infty} }
    { (qt^{-1}q^{z_i-z_j})_{\infty} } \dfrac{
    (tq^{z_i+z_j })_{\infty} } { (q^{z_i+z_j } )_{\infty} } \dfrac{ (t
    q^{z_i -z_j})_{\infty} (q^{1+z_i+z_j} )_{\infty} }
{(q^{z_i-z_j})_{\infty} (qt^{-1} q^{z_i+z_j} )_{\infty} } \right\} \\ \cdot
\sum_{\lambda \subseteq k^n }  
\prod_{i=1}^n (-1)^{\lambda_i}\, q^{4 z_i\lambda_i} t^{(i-n)\lambda_i}
q^{5\lambda_i^2/2 - \lambda_i/2 }   
\dfrac { (t^{n-i}q^{1+\lambda_i})_{\infty} }{
  (q^{2z_i} t^{i-n} q^{\lambda_i})_{\infty} } \dfrac
{(q^{2z_i+2\lambda_i})_{\infty} } {
  (q^{1+2z_i+2\lambda_i})_{\infty} } \\ \cdot (q^{1+k} t^{i-1}
  q^{-\lambda_i })_\infty  (q^{1+k} t^{i-1} q^{2z_i+
    \lambda_i})_\infty  \\ \cdot
  \prod_{1\leq i< j \leq n} \left\{ \dfrac{
      (t^{-1}q^{1+z_i-z_j+\lambda_i - \lambda_j})_{\infty} } {
      (q^{1+z_i-z_j+\lambda_i - \lambda_j})_{\infty} } \dfrac{
      (q^{z_i+z_j + \lambda_i + \lambda_j})_{\infty} } { (tq^{z_i+z_j
        + \lambda_i + \lambda_j})_{\infty} } \right. \\ \left. \cdot \dfrac
{(q^{z_i-z_j+\lambda_i -\lambda_j})_{\infty} (t^{-1}
  q^{1+z_i+z_j+\lambda_i+\lambda_j} )_{\infty} } { (t q^{z_i -z_j+\lambda_i
-\lambda_j} )_{\infty} (q^{1+z_i+z_j+\lambda_i+\lambda_j})_{\infty} }
\right\} 
\end{multline}
and 
\begin{multline}
\label{bsWt2ndversion}
\prod_{1\leq i<j \leq n} \left\{\dfrac{ (q^{1+z_i-z_j})_{\infty} }
  {(qt^{-1}q^{z_i-z_j})_{\infty} } \dfrac{(tq^{z_i-z_j})_{\infty } }
   {(q^{z_i-z_j })_{\infty } } \right\}  
\prod_{i=1}^n  \dfrac{(q^{1+2z_i}t^{i-1})_{\infty} }
{ (q^{1+k} t^{i-1} )_\infty  (qt^{n-i} )_\infty } \\ \cdot
\sum_{\lambda \subseteq k^n } \prod_{i=1}^n
   q^{2z_i\lambda_i} t^{2(i-1)\lambda_i} q^{\lambda_i^2}
  t^{(1-n)\lambda_i} (q^{1+k} t^{i-1} 
  q^{-\lambda_{i}})_\infty (qt^{n-i} q^{\lambda_i})_\infty \\ \cdot
   \prod_{1\leq i<j \leq n}
   \left\{\dfrac{(qt^{-1}q^{z_i-z_j+\lambda_i-\lambda_j} )_\infty} {
   (q^{1+z_i-z_j+\lambda_i-\lambda_j})_\infty } \dfrac {(q^{z_i-z_j +\lambda_i
   -\lambda_j})_\infty } {(tq^{z_i-z_j+\lambda_i -\lambda_j})_\infty }
  \right\}  
\end{multline}
respectively. The Dominated Convergence Theorem will now be applied on
both series above as $k\rightarrow \infty$ to get the
unspecialized \RSi. 

Let $f_\lambda(k)$ denote the summand of the series
in~(\ref{wpsWt2ndversion}) and $g_\lambda(k)$ denote that
in~(\ref{bsWt2ndversion}). Let $h$ be either $f$ or $g$, and consider
the sum   
\begin{equation}
  \sum_{\lambda\in L^+_k} h_\lambda(k)
\end{equation}
where $L^+_k = \{\lambda\in \mathbb{Z}^n: k\geq
\lambda_1 \geq \ldots \lambda_n \geq 0 \}$.
Note that, in the limit, the lattice 
$L^+:=\lim_{k\rightarrow\infty} L^+_k$ consists of all partitions of
length at most $n$.  First verify that the
pointwise limit $ h_\lambda:= \lim_{k\rightarrow\infty}
h_{\lambda}(k)$ exists for all $\lambda\in L^+_k$. Then compute
$m^h_\lambda$ for each $\lambda$ such that $\abs{h_\lambda(k)}
\leq m^h_\lambda$ for any $k$ larger than $\lambda_1$, and finally
verify that the series $\sum_{\lambda\in L^+} m^h_\lambda$ is
convergent. 

That the pointwise limit exists on both sides is clear, since 
\begin{equation}
 \lim_{k\rightarrow \infty} (q^{1+k} t^{i-1}
  q^{-\lambda_i })_\infty  (q^{1+k} t^{i-1} q^{2z_i+
    \lambda_i})_\infty = 1
\end{equation}
on the \wps, and
\begin{equation}
 \lim_{k\rightarrow \infty} (q^{1+k} t^{i-1} 
  q^{-\lambda_{i}})_\infty = 1
\end{equation}
on the \bs, and none of the other factors of $f_\lambda(k)$ or
$g_\lambda(k)$ depend on $k$.

Standard theorems on infinite products and sequences imply that all
the factors of the form $(uq^{\alpha})_\infty /(vq^{\alpha})_\infty$
inside the sum are bounded when $\alpha$ is a non--negative integer
and $u,v\in\mathbb{C}$ such that the denominator never
vanishes, since  
\begin{equation}
  \lim_{\alpha \rightarrow \infty} \left\vert \dfrac{
     (uq^{\alpha})_\infty } {(vq^{\alpha})_\infty } \right\vert = 1
\end{equation}
when $\abs{q}<1$.
Therefore it follows that  
for some constants $C_f$ and $C_g$ that depend only on $q$ and $z$,
and are independent of $k$ and $\lambda$
\begin{equation}
 m^f_\lambda = C_f \abs{\prod_{i=1}^n (-1)^{\lambda_i}\, q^{4
    z_i\lambda_i} t^{(i-n)\lambda_i} 
q^{5\lambda_i^2/2 - \lambda_i/2 } }     
\end{equation}
and 
\begin{equation}
 m^g_\lambda = C_g \abs{\prod_{i=1}^n
   q^{2z_i\lambda_i} t^{2(i-1)\lambda_i} q^{\lambda_i^2}
  t^{(1-n)\lambda_i} } 
\end{equation}
Finally, it needs to be shown that $\sum_{\lambda\in L^+} m^h_\lambda$
is convergent. But this is clearly true for any $t, b\in\mathbb{C}$ 
when $\abs{q}<1$ due to the quadratic factors of $q$. Consider, for
example, the multiple series ratio test. Let $\varepsilon_i =
(0\ldots 1 \ldots 0)$ denote the $n$--tuple of integers with a 1 only
in the $i$--th 
position and zeroes at other positions. Then one sees that 
\begin{equation}
 \abs{ \dfrac{ m^h_{\lambda+\varepsilon_i } } {m^h_\lambda } } 
\end{equation}
is a function of $\abs{q}^{\lambda_i}$ which becomes arbitrarily small
as $\lambda_i \rightarrow \infty$ for each $i\in [n]$ where possible
(i.e., when $\lambda+\varepsilon_i$ is a partition).
Therefore $\sum_{\lambda\in L^+}
m^h_\lambda$ converges when $\abs{q}<1$ as desired. 
\end{proof}

\begin{remark}
The $BC_n$ \RSi~(\ref{unspecializedRSi}) was given in~\cite{Coskun} as
a limiting case of the $BC_n$ \tns\ transformation which was first proved
there (also see~\cite{Coskun1} and~\cite{Rains1}).  
\end{remark}

\subsection{Specializations}
\label{specializations}
Next, it will be shown that the series on both sides of the
\RSi~(\ref{unspecializedRSi}) can be multilateralized, that is they
can be replaced by series over the full lattice
$\mathbb{Z}^n$ under certain specializations of the parameters $b$ and
$t$. An auxiliary result called \ml\ will be needed. 

The root system terminology used in the sequel should be introduced at this
point \cite{Humphreys1}, \cite{Morris1}. 
Let $R$ be the root system $C_n$ in the $n$--dimensional Euclidean
space $E^n$ endowed with the standard inner product $\langle \cdot ,
\cdot \rangle$. Let $R^+$ denote the set of all positive roots, and
$W$ denote the \Wg\ of $R$.
The weight lattice 
\begin{equation}
  L=L(R):=\{\lambda \in E^n : \langle \lambda , \alpha^\vee \rangle
  \in \mathbb{Z}, \forall \, \alpha \in R \}
\end{equation}
and the cone of dominant weights
\begin{equation}
  L^+ = L(R^+):=\{\lambda \in E^n : \langle \lambda , \alpha^\vee
  \rangle \in \mathbb{Z}_{\geq}, \forall \, \alpha \in R^+ \}
\end{equation}
are defined in the standard way. For an arbitrarily fixed
$z\in\mathbb{R}^n$, let $R_z$ denote the isomorphic root system 
obtained by translating $R$ by $-z$ so that the origin is 
moved to $-z$. Let $R^+_z$, $C_z$ 
denote the translated positive roots and fundamental chamber
corresponding to $R^+$ and $C$ of the original root system.  

Now, fix a $z\in\mathbb{R}^n$ such that 
under the reflections $w_{\alpha_z}$ the lattice $L$ remains
fixed. This, of course, implies that 
the standard action of the \Wg\ on the translated root system $R_z$,
by permuting and changing signs of the coordinates, leaves the
lattice $L$ invariant. In other words, it is required that for all $\mu\in
L$ there exists a unique $\lambda\in L$ such that
\begin{equation}
\label{conditionOnz}
  \mu=w(z+\lambda)-z
\end{equation}
for all $w\in W$. This amounts to saying that $w(z)-z \in L$.   

Consider the multiple sum over the dominant cone  
\begin{equation}
\label{SumToMult}
  \sum_{\mu\in L_+} f(z+\mu)
\end{equation}
and assume that the summand has the symmetry $f(\mu+z) =
f(w(\mu+z))$ for all $w\in W$. The interesting cases are when
the sum~(\ref{SumToMult}) is a (convergent) multiple \bhs\ or a
multiple $q$--series. A technique will next be developed showing how
to multilateralize~(\ref{SumToMult}), that is how to replace it, under
certain restrictions, with a series over the entire weight
lattice   
\begin{equation}
  \sum_{\mu\in L} g(z+\mu)
\end{equation}
where $g$ is related to $f$ in a way made precise by 
Lemma~\ref{multilateralization} below.   

Recall that for a reduced root system $R$ Macdonald's polynomial
identity~\cite{Macdonald4} can, in the notation given
above, be written as
\begin{equation}
\label{MacdPolyIden}
  \sum_{w\in W} \prod_{\alpha \in R^+}
   \dfrac{1-u_{\alpha} e^{-w \alpha} }{1- e^{-w\alpha} }  = \sum_{w\in
   W} \prod_{\alpha \in R(w) } u_\alpha  
\end{equation}
where $R(w) = R^+ \cap -wR^+$, $u_\alpha$ are indeterminates
indexed by positive roots $\alpha\in R^+$, and $e^\alpha$ are formal
exponentials denoting elements in the group ring of the root lattice
generated by $R$. The identity will be rewritten in a more convenient
form for the purpose of this paper as
\begin{equation}
\label{genMacdPolyIden}
  \sum_{w\in W} \prod_{\alpha \in R^+}
   \dfrac{(u_{\alpha} e^{-w \alpha} )_\infty } {( e^{-w\alpha})_\infty
   } \dfrac{(q e^{-w \alpha} )_\infty } {( q u_{\alpha}
   e^{-w\alpha})_\infty } = \sum_{w\in W} \prod_{\alpha \in R(w) }
   u_\alpha   
\end{equation}
which follows directly from the original identity~(\ref{MacdPolyIden})
and the definition of the \qPs~(\ref{qPochSymbol}).
Setting $u_\alpha =0$ for all $\alpha\in R^+$ gives a 
rewriting of the \Wdenf\ 
\begin{equation}
\label{WdenfMacIden}
  \sum_{w\in W} \prod_{\alpha \in R^+}
   \dfrac{(q e^{-w \alpha} )_\infty } {( e^{-w\alpha})_\infty } = 1
\end{equation}
Specializing $u_\alpha=1$ for all $\alpha\in R^+$ instead, gives
\begin{equation}
\label{cardWMacIden}
   \lvert\, W \rvert 
= \sum_{w\in W} \prod_{\alpha \in R(w) } u_\alpha  
\end{equation}
It also follows from Macdonald's paper~\cite{Macdonald4}
that the more general $BC_n$ type identity can be written as 
\begin{equation}
\label{genMacdPolyIdenBCn}
  \sum_{w\in W} \prod_{\alpha \in R^+}
   \dfrac{(u^{1/2}_{2\alpha} u_{\alpha} e^{-w \alpha} )_\infty }
   {(u^{1/2}_{2\alpha} e^{-w\alpha})_\infty 
   } \dfrac{(q u^{1/2}_{2\alpha} e^{-w \alpha} )_\infty } {( q
   u^{1/2}_{2\alpha} u_{\alpha}
   e^{-w\alpha})_\infty } = \sum_{w\in W} \prod_{\alpha \in R(w) }
   u_\alpha   
\end{equation}
where $u^{1/2}_{2\alpha}=1$ by convention, when $2\alpha$ is not a
root. If one sets 
$u^{1/2}_{2\alpha}=-1$ for all $\alpha\in R_s^+$, the short roots of
$R^+$, the \lhs\ of the identity~(\ref{genMacdPolyIdenBCn}) becomes   
\begin{equation}
\label{genMacdPolyIdenBn-}
  \sum_{w\in W} \prod_{\alpha \in R_m^+}
   \dfrac{(u_{\alpha} e^{-w \alpha} )_\infty }
   {( e^{-w\alpha})_\infty } \dfrac{(q e^{-w \alpha} )_\infty } {( q
   u_{\alpha} e^{-w\alpha})_\infty } \prod_{\alpha \in R_s^+}
   \dfrac{(-u_{\alpha} e^{-w \alpha} )_\infty }
   {(- e^{-w\alpha})_\infty } \dfrac{(-q e^{-w \alpha} )_\infty } {( -q
   u_{\alpha} e^{-w\alpha})_\infty } 
\end{equation}
where $R_m^+$ denotes the medium positive roots. 

Some more terminology will be needed in the proof below. Under the
action of the \Wg\ on $R_z$, the image of $L^+$ may
not be the full weight lattice $L$. Depending on $z$, some subsets of
$L$ may be mapped more than once, and some subsets may not be the
image of any subset of $L^+$ under this action. 

Let $P$, called an overlap, denote the (possibly empty) subset of $L$
defined by 
\begin{equation}
  \label{def:overlap}
  P = \{\lambda\in L : \lambda\in w_i L^+ \cap w_j L^+,
  \mathrm{\,\,where\,\,} i\neq j \mathrm{\,\,and\,\,} i,j\in [2^n n!]\}   
\end{equation} 
In order to avoid pathological situations with complicated overlap
sets, a further restriction is imposed on $z\in \mathbb{R}^n$, and
it is required that the $-z$ 
(the new center) is located outside the fundamental chamber $C$ so
that the translated fundamental chamber $C_z$ properly contains
$C$. Under this restriction, the overlap $P$ is non--empty only 
when $-z$ is on one of the walls of $C$, and 
the overlap $P$ turns out to be a proper subset of the walls of the
Weyl chambers for the root system $R$.

Let $Q$, called a gap, denote the (possibly empty) subset of $L$
containing all $\mu\in L$ that is not an image of any $\lambda\in L^+$
under the standard action of the \Wg\ on $R_z$. In other
words, $Q$ is the set theoretical difference   
\begin{equation}
  \label{def:gap}
  Q = L \backslash \cup_{w\in W} wL^+ 
\end{equation} 

The technique that will be used in multilateralizing \bhs\ or
$q$--series associated to root systems is developed in the next lemma. 
\begin{lemma}[Multilateralization]
\label{multilateralization}
Suppose that the summand $f(\zeta+\mu)$ of the sum 
\begin{equation}
  \sum_{\mu\in L_+} f(\zeta+\mu)
\end{equation}
has the invariance property that $f(\zeta+\mu) =f(w(\zeta+\mu))$ for
any $\zeta \in\mathbb{R}^n$. Fix $z\in\mathbb{R}^n$ 
such that $L$ is invariant under the action of the \Wg\ on $R_z$ 
as defined by condition~(\ref{conditionOnz}). Choose the fundamental Weyl
chamber $C$, so that $L^+\subseteq C$ and further require that 
$C\subseteq C_z$. Another obvious restriction on the choice of $z$ is
that $f(z+\mu)$ is well defined (has no essential poles) at any
$\mu\in L^+$. 
Let $P$ denote the overlap and $Q$ denote the gap \wrt\ this action.
Suppose further that
\begin{equation}
\label{gap}
  \sum_{\mu\in Q} g(z+\mu) =0 
\end{equation}
and
\begin{equation}
\label{overlap}
  \bigg(\sum_{w\in W} \prod_{\alpha \in R(w) } u_\alpha \bigg)
  \sum_{\mu\in P\cap L^+} f(z+\mu) =  \sum_{\mu\in P} g(z+\mu)  
\end{equation}
where $g$ is defined by 
\begin{equation}
   g(z+\mu) := f(z+\mu) \bigg( \prod_{\alpha \in R^+}
   \dfrac{(u^{1/2}_{2\alpha} u_{\alpha} q^{\langle \alpha, \,z+\mu
   \rangle } )_\infty } {(u^{1/2}_{2\alpha} q^{\langle \alpha, \, z+\mu
   \rangle } )_\infty 
   } \dfrac{(q u^{1/2}_{2\alpha} q^{\langle \alpha, \, z+\mu
   \rangle } )_\infty } {( q u^{1/2}_{2\alpha} u_{\alpha}
   q^{\langle \alpha, \, z+\mu \rangle } )_\infty } \bigg)    
\end{equation}
in terms of the indeterminates $u_\alpha$. 
Under these conditions, 
\begin{equation}
  \bigg(\sum_{w\in W} \prod_{\alpha \in R(w) } u_\alpha \bigg)
  \sum_{\mu\in L^+} f(z+\mu) = \sum_{\mu\in L}\, g(z+\mu)  
\end{equation}
\end{lemma}

\begin{proof}
The identity~(\ref{genMacdPolyIdenBCn}) can be written in the form
\begin{equation}
\label{genMacdPolyIdenBCn2ndversion}
  \sum_{w\in W} \prod_{\alpha \in R^+}
   \dfrac{(u^{1/2}_{2\alpha} u_{\alpha} q^{- \langle \alpha, \,w(z+\mu)
   \rangle } )_\infty } {(u^{1/2}_{2\alpha} q^{-\langle \alpha, \,
   w(z+\mu) \rangle } )_\infty 
   } \dfrac{(q u^{1/2}_{2\alpha} q^{-\langle \alpha, \, w(z+\mu)
   \rangle } )_\infty } {( q u^{1/2}_{2\alpha} u_{\alpha}
   q^{-\langle \alpha, \, w(z+\mu) \rangle } )_\infty }
= \sum_{w\in W} \prod_{\alpha \in R(w) } u_\alpha   
\end{equation}
Next, write 
\begin{equation}
\label{separateOverlap}
  \sum_{\mu\in L^+} f(z+\mu) = \sum_{\mu\in P\cap L^+} f(z+\mu) + 
  \sum_{\mu\in L^+\backslash P} f(z+\mu) 
\end{equation}
and consider the second sum. Multiplying $\sum_{\mu\in
  L^+\backslash P } f(z+\mu)$ by the \rhs\ of the above 
identity~(\ref{genMacdPolyIdenBCn2ndversion}) one gets
\begin{multline}
  \bigg( \sum_{w\in W} \prod_{\alpha \in R(w) } u_\alpha \bigg)
  \sum_{\mu\in L^+\backslash P } 
  f(z+\mu) \\ 
 =  \sum_{\mu\in L^+\backslash P} \sum_{w\in W} f(w(z+\mu))
  \prod_{\alpha \in R^+}
   \dfrac{(u^{1/2}_{2\alpha} u_{\alpha} q^{- \langle \alpha, \,w(z+\mu)
   \rangle } )_\infty } {(u^{1/2}_{2\alpha} q^{-\langle \alpha, \,
   w(z+\mu) \rangle } )_\infty 
   } \dfrac{(q u^{1/2}_{2\alpha} q^{-\langle \alpha, \, w(z+\mu)
   \rangle } )_\infty } {( q u^{1/2}_{2\alpha} u_{\alpha}
   q^{-\langle \alpha, \, w(z+\mu) \rangle } )_\infty }
\end{multline}
Note that $f(z+\mu)$ is replaced by $f(w(z+\mu))$ on the \rhs\ due to
the $W$ invariance of the summand $f$. 
If the order of summation is switched and the sum is written as
$\sum_{w\in W} 
\sum_{\mu\in L^+\backslash P}$, it becomes clear from the definitions
that this double sum can be written as a single sum over the subset
$L\backslash ( Q \cup P ) $. Because, this is precisely the set that
contains all $\lambda\in L$ such that $w(\mu+z) =\lambda+z$ for some
$\mu \in L^+\backslash P$ and~$w\in W$. Thus, it follows that 
\begin{multline}
  \bigg( \sum_{w\in W} \prod_{\alpha \in R(w) } u_\alpha \bigg)
  \sum_{\mu\in L^+\backslash P } 
  f(z+\mu) \\ 
 =  \sum_{\mu\in L\backslash ( Q \cup P ) } f(z+\mu)  \bigg(
  \prod_{\alpha \in R^+}
   \dfrac{(u^{1/2}_{2\alpha} u_{\alpha} q^{- \langle \alpha, \, z+\mu
   \rangle } )_\infty } {(u^{1/2}_{2\alpha} q^{-\langle \alpha, \,
   z+\mu \rangle } )_\infty 
   } \dfrac{(q u^{1/2}_{2\alpha} q^{-\langle \alpha, \, z+\mu
   \rangle } )_\infty } {( q u^{1/2}_{2\alpha} u_{\alpha}
   q^{-\langle \alpha, \, z+\mu \rangle}  )_\infty } \bigg)    
\end{multline}
But due to the condition~(\ref{gap}) this sum on the \rhs\ can actually
be written over $L\backslash P$. Therefore, assuming that the
condition~(\ref{gap}) holds, one gets 
\begin{equation}
  \bigg( \sum_{w\in W} \prod_{\alpha \in R(w) } u_\alpha \bigg)
 \sum_{\mu\in L^+\backslash P } f(z+\mu)  
 =  \sum_{\mu\in L\backslash P } g(z+\mu) 
\end{equation}
Finally adding the sum $\sum_{\mu\in P} g(z+\mu) $ to both sides and 
using the condition~(\ref{overlap}) gives
\begin{equation}
  \bigg( \sum_{w\in W} \prod_{\alpha \in R(w) } u_\alpha \bigg)
  \sum_{\mu\in L^+ } f(z+\mu) =  \sum_{\mu\in L } g(z+\mu)  
\end{equation}
as desired.
\end{proof}

Note that the argument used in the proof of Lemma~\ref{multilateralization}
generalizes the technique used in Section \ref{AGeneralRRis} for the 
alternative proof of a generalization of the \RRis. 

In order to show that the multiple \RSi~(\ref{unspecializedRSi}) can be
multilateralized, the symmetries of the summands in both sides
should be verified. 
\begin{lemma}
\label{symmetriesOfRSi}
The summand of the series in the \wps\ of the
\RSi~(\ref{unspecializedRSi}) has a $C_n$ type symmetry, 
and the one in the balanced side has an $A_n$ type symmetry. 
\end{lemma}

\begin{proof}
That the series in the \wps\ of~(\ref{unspecializedRSi}) has $C_n$
type symmetries is already established in a more general setting in
Remark~\ref{HyperoctahedralSymmetry}. 

The series in the \bs\ is invariant under the standard action of the
Weyl group of the root system $A_n$ of rank $n$ (i.e., the symmetric
group in $n$ letters). This can be seen similarly by setting $q^{z_i}
= b^{1/2} t^{1-i}$ and again using (\ref{qPochSymbol}),
(\ref{qPsFlipIdentity}) and  
\begin{equation}
\label{doubleProductTwo}
\prod_{1\leq i<j \leq n} a_i 
=  \prod_{ i=1 }^{n} a_i^{n-i} 
\end{equation}
to rewrite the series in the following form. 
\begin{multline}
 \prod_{i=1}^n  q^{-z_i^2 }
   \dfrac{(q^{1+2z_i}t^{i-1})_{\infty} } { 
   (qt^{n-i})_{\infty} }  \prod_{1\leq i<j \leq n} 
   \left\{ \dfrac{(q^{1+z_i-z_j})_{\infty} }
   {(t^{-1}q^{1+z_i-z_j})_{\infty} } 
  \dfrac{(q^{1+z_j-z_i})_{\infty} } {(t^{-1}q^{1+z_j-z_i})_{\infty} }
  \right\}  \\  \cdot  \sum_{\mu,\ell(\mu)\leq n} \prod_{i=1}^n
   \left\{ q^{(z_i+\mu_i)^2 } (b^{-1/2} t^{n-1} q^{1+\mu_i+z_i}
   )_\infty \right\} \\  \cdot 
 \prod_{1\leq i<j \leq n}
   \left\{ \dfrac {(t^{-1}q^{1+z_i+\mu_i-z_j-\mu_j})_{\infty}}
   {(q^{1+z_i+\mu_i-z_j-\mu_j} )_{\infty} }
  \dfrac {(t^{-1}q^{1-z_i-\mu_i+z_j+\mu_j})_{\infty} }
   {(q^{1-z_i-\mu_i+z_j+\mu_j} )_{\infty} }  \right\}
\end{multline}
It is obvious from this representation that the series on the \bs\ is
invariant under the permutations $q^{z+\mu} \leftrightarrow
q^{w(z+\mu)}$ of the symmetric group.  
\end{proof}

\begin{remark}
\label{centerz}
The crucial point is to fix a center $-z\in\mathbb{Z}^n$ for the
\ml\ satisfying all the restrictions imposed by 
Lemma~\ref{multilateralization}. The
condition~(\ref{conditionOnz}) implies $w(z)-z \in L$, which
in turn implies that 
each $z_i$ is of the form $m_i/2$ for integers $m_i$ and that all
$m_i$ for $i\in [n]$ have the same parity. This observation together
with the hypothesis $C\subseteq
C_z$ in Lemma~\ref{multilateralization} implies that
$z=m/2+\mu=(m/2+\mu_1,\ldots, m/2+\mu_n)$ where
$m$ is a non--negative integer and $\mu$ is some partition with at
most $n$ parts and $\mu_n=0$. On the other hand, the substitution 
$q^{z_i} = b^{1/2} t^{1-i}$ made in the 
series in~(\ref{unspecializedRSi}) yields that $z_i=
1/2\log_q b + (1-i) \log_q t$. 

Combine the two equivalent representations  $z_i=
1/2\log_q b + (1-i) \log_q t$ and $z_i=m/2+\mu_i$ from the last
paragraph to show that  
multilateralization is possible when $t=q^k$ for some non--negative
integer $k=\mu_{n-1}$. It also turns out that $b=q^{m+2(n-1)k}$ and,
in general, $z$ has to have the form  
\begin{equation}
  \label{eq:specializations}
  z_i = m/2+k(n-i)
\end{equation}
where $m$ and $k$ are non--negative integers.  

It is worthwhile to note that the assumption 
$C\subseteq C_z$ is not essential in the proof of
Lemma~\ref{multilateralization} 
above. It guarentees that the overlap and gap sets are not
too complicated in applications. In principle, this assumption may be
dropped leading, in particular,  
to specializations where $m$ and $k$ could be negative. In fact, an
example is presented in Section~\ref{AGeneralRRis} where $m$ was an
arbitrary integer.  
\end{remark}

The following notation will be used in the sequel. $\aleph_\lambda$ is
defined by 
\begin{equation}
  \label{def:delta}
  \aleph_\lambda = \begin{cases} 1, \mathrm{\;if \; \lambda \;is \;a\;
  rectangular\; partition} \\ n!, \mathrm{otherwise} \end{cases} 
\end{equation}
and the product $\sideset{}{^{r,s}} \prod$ denotes
\begin{equation}
  \label{def:prod}
  \sideset{}{^{r,s}} \prod f(x) = \begin{cases} s, \mathrm{\,if \,}
  r=0 \\ \prod f(x), \mathrm{otherwise} \end{cases}    
\end{equation}

\begin{lemma}[$B_n$ and $D_n$ Specializations]
\label{specializationsOfRSi}
Set $q^{z_i} = b^{1/2} t^{1-i}$ for $i\in [n]$ in the unspecialized
\RSi~(\ref{unspecializedRSi}). With the specialization $z_i=
m/2+k(n-i)$, where $m$ and $k$ are non--negative
integers,
the series on the \rhs\ of the identity~(\ref{unspecializedRSi}) can
be written as follows. 

\begin{multline}
\label{specializedRSi}
\prod_{i=1}^n (q^{1+k(n-i)})_{\infty} 
 \cdot \sum_{\lambda, \ell(\lambda)<n }  \prod_{i=1}^{n} \left\{ 
   q^{(m+k(n-1) ) \lambda_i + \lambda_i^2}
\dfrac{1 } {(q^{1+k(n-i)})_{\lambda_i }} \right\} \\  \cdot  
  \sideset{}{^{k,\aleph_\lambda}} \prod_{1\leq i<j \leq n}
   \left\{\dfrac{(q^{1+k(-1-i+j)+\lambda_i-\lambda_j} )_\infty} {
   (q^{1+k(-i+j)+\lambda_i-\lambda_j})_\infty } \dfrac{(q^{k(-i+j)
   +\lambda_i -\lambda_j})_\infty } {(q^{k(1-i+j) +\lambda_i
   -\lambda_j})_\infty } \right\}   \\
=  \sum_{\lambda \in \mathbb{Z}^n }  
\prod_{i=1}^n \left\{ (-1)^{\lambda_i}\, q^{(-1/2+2m+3k(n-i))\lambda_i +
   5\lambda_i^2/2 } \right\}
   \sideset{}{^{m,1}}\prod_{i=1 }^n  
\left\{ \dfrac{(q^{1+k(n-i) +\lambda_i})_{\infty} }{  
  (q^{m+ k(n-i) + \lambda_i})_{\infty} } \right\} \\  \cdot
  \sideset{}{^{k,1}}\prod_{1\leq i< j \leq n} \left\{
   \dfrac{ (q^{1+ k(-1-i+j) + 
   \lambda_i - \lambda_j})_{\infty} } {(q^{m+ k(1+2n-i -j ) 
        + \lambda_i + \lambda_j})_{\infty} } 
    \dfrac{ ( q^{1+ m + k(-1+2n-i-j) +\lambda_i+\lambda_j} )_{\infty} } {
      (q^{ k(1-i+j) +\lambda_i -\lambda_j} )_{\infty} }
\right\} 
\end{multline}
The $m>1$ case corresponds to the $B_n$ specializations and that of
$m\in\{0,1\}$ corresponds to the $D_n$ specializations. 
\end{lemma}

\begin{proof}
It follows from~(\ref{wpsWt2ndversion}) and~(\ref{bsWt2ndversion}) that
the \bs\ and the \wps\ of the \RSi~(\ref{unspecializedRSi}) can be
written in the form 
\begin{multline}
\label{RSi2ndversion}
\prod_{i=1}^n (q^{1+2z_i}t^{i-1})_{\infty} \cdot 
\sum_{\lambda,\ell(\lambda)\leq n} \prod_{i=1}
   q^{2z_i\lambda_i+\lambda_i^2} t^{(-1-n+2i)\lambda_i} 
    \dfrac{1 } {(qt^{n-i})_{\lambda_i}} \\  \cdot  
   \prod_{1\leq i<j \leq n}
   \left\{\dfrac{(qt^{-1}q^{z_i-z_j+\lambda_i-\lambda_j} )_\infty} {
   (q^{1+z_i-z_j+\lambda_i-\lambda_j})_\infty } \dfrac {(q^{z_i-z_j +\lambda_i
   -\lambda_j})_\infty } {(tq^{z_i-z_j+\lambda_i -\lambda_j})_\infty }
\right\}   \\
= \prod_{i=1}^n \dfrac{ (q^{2z_i}t^{i-n})_{\infty} } {
  (qt^{n-i})_{\infty} } 
\dfrac{(q^{1+2z_i} )_{\infty}} { (q^{2z_i})_{\infty} } \cdot  
\prod_{1\leq i< j \leq n} \left\{ \dfrac{
    (tq^{z_i+z_j })_{\infty} } { (q^{z_i+z_j } )_{\infty} } 
\dfrac{ (q^{1+z_i+z_j} )_{\infty} }
{ (qt^{-1} q^{z_i+z_j} )_{\infty} } \right\} \\ 
 \cdot  \sum_{\ell(\lambda)\leq n}  
\prod_{i=1}^n (-1)^{\lambda_i}\, q^{4 z_i\lambda_i} t^{(i-n)\lambda_i}
q^{5\lambda_i^2/2 - \lambda_i/2 }   
\dfrac { (t^{n-i}q^{1+\lambda_i})_{\infty} }{
  (q^{2z_i} t^{i-n} q^{\lambda_i})_{\infty} } \dfrac
{(q^{2z_i+2\lambda_i})_{\infty} } {
  (q^{1+2z_i+2\lambda_i})_{\infty} }  \\  \cdot 
  \prod_{1\leq i< j \leq n} \left\{ \dfrac{
      (t^{-1}q^{1+z_i-z_j+\lambda_i - \lambda_j})_{\infty} } {
      (q^{1+z_i-z_j+\lambda_i - \lambda_j})_{\infty} } \dfrac{
      (q^{z_i+z_j + \lambda_i + \lambda_j})_{\infty} } { (tq^{z_i+z_j
        + \lambda_i + \lambda_j})_{\infty} } \right. \\ \left.  \cdot
    \dfrac{(q^{z_i-z_j+\lambda_i -\lambda_j})_{\infty} (t^{-1}
  q^{1+z_i+z_j+\lambda_i+\lambda_j} )_{\infty} } { (t q^{z_i -z_j+\lambda_i
-\lambda_j} )_{\infty} (q^{1+z_i+z_j+\lambda_i+\lambda_j})_{\infty} }
\right\} 
\end{multline}
The series on the \rhs\ of~(\ref{RSi2ndversion}) will be
multilateralized. Remark~\ref{centerz} preceding this lemma
explain, in the notation of 
Lemma~\ref{multilateralization}, that $L=\mathbb{Z}^n$ is
invariant under the standard action of the \hgr\ $W$ (the semidirect
product of the symmetric group $\mathbb{S}_n$ and $\mathbb{Z}^n_2$) on
$R_z$. That is, under the maps 
$q^{z_i+\mu_i} \leftrightarrow q^{\pm (z_j + \mu_j)}$ for all $i,j\in [n]$,
the lattice $L=\mathbb{Z}^n$ remains invariant when
$z_i=m/2+k(n-i)$ and $t=q^k$.

First note that the front factors on both sides, under this
specialization, combine to give  
\begin{multline}
\label{front}
\prod_{i=1}^n
   \dfrac{1}{(q^{1+k(n-i)})_{\infty}} =  \!\!\! \prod_{1\leq i< j \leq
   n} \!\!\! \left\{ \dfrac{ 
    (q^{m+k(1+2n-i-j) })_{\infty} } { (q^{m+ k(2n-i-j) } )_{\infty} } 
\dfrac{ (q^{1+ m+ k(2n-i -j)} )_{\infty} }
{ (q^{1+m+ k(-1+2n-i -j)} )_{\infty} } \right\} \\
\cdot \prod_{i=1}^n \left\{ \dfrac{
   (q^{m+k(n-i)} )_{\infty} } { (q^{1+k(n-i)})_{\infty} }  
\dfrac{(q^{1+ m+2k(n-i)} )_{\infty}} { (q^{m+2k(n-i)})_{\infty} }
\dfrac{1}{ (q^{1+m+k(2n-i-1)} )_{\infty}} \right\} 
\end{multline}
for any positive integer $m$.

Since the specializations with 
either $m=0$ or $k=0$ produce non--empty overlap
sets, the proof will be
divided into four cases depending on whether $m$ or $k$ is zero. For
each case the subsets $P$ and $Q$ of $L$ and corresponding 
Macdonald polynomial identity will be identified, and the conditions
of Lemma~\ref{multilateralization} will be verified. 

Two types of subsets of the weight lattice $L=\mathbb{Z}^n$ are
special in the application of Lemma~\ref{multilateralization}. They
are defined as the intersection of $L$ and the following two finite
affine hyperplane arrangements:
\begin{enumerate}
\item\label{type1Hyp} $x_i = -m - k(n-i) +1 + \ell$ where $i\in[n]$,
   and $\ell\geq -1$ is an integer. 
\item\label{type2Hyp}
\begin{enumerate}
\item $x_i = x_j - k(j-i+1) + 1 +\ell$, and
\item $x_i = -x_j - m -k(2n-j-i+1) +1 +\ell$
\end{enumerate}
\noindent
where $1\leq i<j\leq n$, and $\ell\geq -1$ is an integer.  
\end{enumerate}

\emph{Case 1: $m,k\neq 0$.} In this case the gap $Q$ consists of the
Type~\ref{type1Hyp} subsets for all $\ell = 0, 1, \ldots, m-2$ and the
Type~\ref{type2Hyp} subsets for all $\ell=0, 1, \ldots, 2(k-1)$. 
The overlap $P$ is empty, that is $P=\emptyset$.   

\emph{Case 2: $m=0, k\neq 0$.} Here, the gap $Q$ consists of the
Type~\ref{type2Hyp} subsets for all $\ell=0, 1, \ldots, 2(k-1)$. The
overlap $P$ is the set theoretical difference of the
Type~\ref{type1Hyp} subsets with $\ell=-1$ and 
the gap $Q$. 

\emph{Case 3: $m\neq 0, k=0$.} Similarly, in this case the gap $Q$
consists of the Type~\ref{type1Hyp} subsets for all $\ell = 0, 1, \ldots,
m-2$. The overlap $P$ is the set
theoretical difference of the Type~\ref{type2Hyp} subsets with
$\ell=-1$ and the gap $Q$.

\emph{Case 4: $m=0, k=0$.} Finally, the gap $Q$ is empty if both $m$
and $k$ are zeros. The overlap $P$ consists of the Type~\ref{type1Hyp}
and Type~\ref{type2Hyp} subsets with $\ell=-1$.

For the cases where $m\neq 0$, the $C_n$ type Macdonald's polynomial
identity~(\ref{genMacdPolyIden}) will be used with $u_\alpha=0$ for
all positive roots $\alpha$. This identity takes on the form
\begin{equation}
\label{specMacdPolyIden}
 \sum_{w\in W} \prod_{1\leq i < j \leq n} \left\{ 
   \dfrac{(q w(x_i^{-1} x_j) )_\infty } {( w(x_i^{-1} x_j ) )_\infty }
   \dfrac{(q w(x_i^{-1} x_j^{-1} ) )_\infty } {(w(x_i^{-1} x_j^{-1} )
   )_\infty } \right \} \prod_{i=1}^n \dfrac{(q w(x_i^{-2}) )_\infty }
   {( w(x_i^{-2} ) )_\infty }  =1
\end{equation}
where $x_i=q^{z_i+\lambda_i}=q^{m/2+k(n-i)+\lambda_i}$ and $W$ is the
\hgr. Setting $u_\alpha=0$ 
will remove a set of factors corresponding to the positive roots for $C_n$
and therefore yields the $B_n$ type specialization. 

When $m=0$, however, the $BC_n$ type 
identity~(\ref{genMacdPolyIdenBn-}) will be used with $u_\alpha=0$ for
all roots of medium length. The identity can be written as
\begin{equation}
\label{specMacdPolyIdenBn-}
 \sum_{w\in W} \prod_{1\leq i < j \leq n} \left\{ 
   \dfrac{(q w(x_i^{-1} x_j) )_\infty } {( w(x_i^{-1} x_j) )_\infty }
   \dfrac{(q w(x_i^{-1} x_j^{-1} ) )_\infty } {(w(x_i^{-1} x_j^{-1} )
   )_\infty } \right \} \prod_{i=1}^n \dfrac{(-q w(x_i^{-1}) )_\infty }
   {( -w(x_i^{-1} ) )_\infty }  =1
\end{equation}
where $x_i$ and $W$ are defined as before. 
This choice will remove a set of factors corresponding to the long,
the short and the medium positive roots for $BC_n$ 
and therefore yields the $D_n$ type specialization. 

In both cases, one gets
\begin{multline}
g(z+\lambda) \!\! = 
\!\!\!  {\sideset{}{^{k,1 }} \prod_{1\leq i< j \leq n}} \!\! \left\{\!\!
   \dfrac{ (q^{1+ k(-1-i+j) +  \lambda_i - \lambda_j})_{\infty} }
   {(q^{m+ k(1+2n-i -j )  + \lambda_i + \lambda_j})_{\infty} } 
    \dfrac{ ( q^{1+ m + k(-1+2n-i-j) +\lambda_i+\lambda_j} )_{\infty} } {
      (q^{ k(1-i+j) +\lambda_i -\lambda_j} )_{\infty} }
\!\! \right\} \\
\cdot \prod_{i=1}^n \left\{ (-1)^{\lambda_i}\,
   q^{(-1/2+2m+3k(n-i))\lambda_i + 
   5\lambda_i^2/2 } \right\} \sideset{}{^{m,1}}\prod_{i=1}^{n} \left\{
   \dfrac{ (q^{1+k(n-i) +\lambda_i})_{\infty} }{ 
  (q^{m+ k(n-i) + \lambda_i})_{\infty} } \right\} 
\end{multline}
where the summand $g(z+\lambda)$ is as defined in
Lemma~\ref{multilateralization}.  

The rest of the proof shows that the gap condition~(\ref{gap})
and the overlap condition~(\ref{overlap}) in
Lemma~\ref{multilateralization} are both satisfied. The factors of the
form 
\begin{equation}
\sideset{}{^{m,1}}\prod_{i=1}^n \left\{ \dfrac{ (q^{1+k(n-i)
  +\lambda_i})_{\infty} }{ (q^{m+ k(n-i) + \lambda_i})_{\infty} }
  \right\} \mathrm{,} \quad
\sideset{}{^{k, 1}}\prod_{1\leq i< j \leq n} \left\{ 
    \dfrac{ ( q^{1+ m + k(-1+2n-i-j) +\lambda_i+\lambda_j} )_{\infty}
  } {(q^{m+ k(1+2n-i -j ) + \lambda_i + \lambda_j})_{\infty} }  
\right\} 
\end{equation}
and
\begin{equation}
\sideset{}{^{k, 1} }\prod_{1\leq i< j \leq n} \left\{
   \dfrac{ (q^{1+ k(-1-i+j) + \lambda_i - \lambda_j})_{\infty} } {
      (q^{ k(1-i+j) +\lambda_i -\lambda_j} )_{\infty} }
\right\} 
\end{equation}
in~(\ref{RSi2ndversion}) show that the summand $g(z+\lambda)$ vanishes
over Type 1, Type 2 (a) and Type 2 (b) subsets, respectively,
verifying the gap condition~(\ref{gap}) in all four cases where the
gap $Q$ is non--empty.     

In all three cases where the overlap $P$ is non--empty, 
the series over
$P\cap L^+$ is removed from the series over $L^+$, and 
the resulting series 
\begin{equation}
\label{twoseries}
  \sum_{\mu\in L^+\backslash (P\cap L^+)} f(z+\mu) \quad \mathrm{and}
  \quad \sum_{\mu\in P\cap L^+} f(z+\mu)
\end{equation}
are multilateralized separately using
Lemma~\ref{multilateralization}, giving 
\begin{equation}
\label{relationfor2nd}
  \sum_{\mu\in L^+\backslash (P\cap L^+)} f(z+\mu) = \sum_{\mu\in
  L\backslash P} g(z+\mu)
\end{equation}
and 
\begin{equation}
\label{relationfor1st}
  \sum_{\mu\in P\cap L^+} f(z+\mu) =  \sum_{\mu\in P} g(z+\mu)  
\end{equation}
For the second series over $P\cap L^+$ 
in~(\ref{twoseries}), a special 
version of Lemma~\ref{multilateralization} is needed. Namely, $L$ is 
replaced by $P$ and 
$L^+$ is replaced by $P\cap L^+$ in the lemma, and appropriate
versions of Macdonald's identities~(\ref{specMacdPolyIden})
and~(\ref{specMacdPolyIdenBn-}) are used depending on whether $m=0$
or $m> 0$ by specializing $u_\alpha$ parameters accordingly.
Finally, the two multilateral series are added together to get the
single multilateral series $\sum_{\mu\in L} g(z+\mu)$.
\end{proof}

\begin{remark}
\label{bothsides_specializationsOfRSi}
The main argument in the proof of of Lemma~\ref{multilateralization}
can be applied to the 
series on the \bs\ of \RSi~(\ref{unspecializedRSi}) to write it as a
sum over $\mathbb{Z}_{\geq}^n$. Here, one uses the following
special case of the $A_n$ version of Macdonald's polynomial
identity~(\ref{genMacdPolyIden})
\begin{equation}
 \sum_{w\in W} \prod_{1\leq i < j \leq n} 
   \dfrac{(q w(x_i^{-1} x_j) )_\infty } {( w(x_i^{-1} x_j^{-1})
   )_\infty } = 1
\end{equation}
where $W=S_n$, the symmetric group on $n$ letters. This is because the
series on the \bs\ has $A_n$ type symmetries. Multiplying both sides
by $\prod_{i=1}^n (q)_{k(n-i)}$ using 
\begin{equation}
  (q)_{\lambda_i} = (q)_{k(n-i)} (q^{1+k(n-i)} )_{\lambda_i-k(n-i) }
\end{equation}
the specialized $BC_n$ \RSi, therefore, can also be written as  
\begin{multline}
\label{specializedRSi_bothsides}
\sum_{\lambda \in \mathbb{Z}_{\geq}^n } \prod_{i=1}^{n} \left\{
\dfrac{ q^{(m+k(n-1) ) (\lambda_i-k(n-i)) + (\lambda_i-k(n-i))^2 } }
{(q)_{\lambda_i } } \right\} \; 
   \sideset{}{^{k, 1}}\prod_{1\leq i<j \leq n}
   \left\{\dfrac{(q^{1-k+\lambda_i-\lambda_j} )_\infty} {
   (q^{k +\lambda_i -\lambda_j})_\infty } \right\}   \\
= \dfrac{1}{(q)^n_{\infty} } \cdot \sum_{\lambda \in \mathbb{Z}^n }  
\prod_{i=1}^n \left\{ (-1)^{\lambda_i}\, q^{(-1/2+2m+3k(n-i))\lambda_i +
   5\lambda_i^2/2 } \right\}
   \sideset{}{^{m,1}}\prod_{i=1 }^n  
\left\{ \dfrac{(q^{1+k(n-i) +\lambda_i})_{\infty} }{  
  (q^{m+ k(n-i) + \lambda_i})_{\infty} } \right\} \\  \cdot 
  \sideset{}{^{k, 1}} \prod_{1\leq i< j \leq n} \left\{
    \dfrac{ (q^{1+ k(-1-i+j) + \lambda_i - \lambda_j})_{\infty} }
    {(q^{m+ k(1+2n-i -j ) + \lambda_i + \lambda_j})_{\infty} } 
    \dfrac{ ( q^{1+ m + k(-1+2n-i-j) +\lambda_i+\lambda_j} )_{\infty} } {
      (q^{ k(1-i+j) +\lambda_i -\lambda_j} )_{\infty} }
\right\} 
\end{multline}
\end{remark}

\subsection{Multiple Rogers--Ramanujan Identities}
\label{RRidentities}

Recall from Section~\ref{AGeneralRRis} that the two bilateralizations
obtained for the \od\ \RSi\  
yielded the two \RRis~(\ref{eq:RRs}) 
by an application of the \Jtpi~(\ref{eq:JtpiBackground}).
A well--known generalization of the \Jtpi\ to arbitrary root
systems is called \Mis~\cite{Macdonald3},~\cite{Gustafson2}.
One expects that a similar procedure would produce product
representations for the specialized $B_n$ and $D_n$
\RSis~(\ref{specializedRSi}) or~(\ref{specializedRSi_bothsides})
using the more general \Mis. This proves to be a non--trivial
problem except in a special case when $k=0$. 
It is possible, however, to give interesting generalizations of  
the classical \RRis\ in terms of determinants of theta functions for
the root systems $B_n$ and $D_n$ of rank $n$, at least for the case
$k=1$.  

First, it will be noted that the trivial specializations corresponding
to $k=0$ (i.e., $t=1$) result in product representations. 
\begin{theorem} 
\label{thm:trivialRRis}
The $D_n$ specializations~(\ref{specializedRSi}) of the \RSi\ give
$n$-fold products of the first and second \RRis~(\ref{eq:RRs}),
respectively. Namely,
\begin{equation}
\label{trivialcase}
\sum_{\lambda \in \mathbb{Z}_{\geq}^n } \prod_{i=1}^{n} \left\{
\dfrac{ q^{\delta \lambda_i  + \lambda_i^2 } }
{(q )_{\lambda_i }} \right\} 
= \prod_{i=1 }^n \left\{ 
  \dfrac{1}{(q^{1+\delta}; q^5)_\infty (q^{4-\delta}; q^5)_\infty}
  \right\} 
\end{equation}
where $\delta\in\{0,1\}$ and $\abs{q}<1$ as usual.
\end{theorem}

\begin{proof}
The proof follows, immediately, from the $k=0$ case
of~(\ref{specializedRSi}) 
and~(\ref{specializedRSi_bothsides}) which can be written as
\begin{multline}
\sum_{\lambda,\; \ell(\lambda)\leq n } \aleph_\lambda \cdot
\prod_{i=1}^{n} \left\{  q^{m \lambda_i + \lambda_i^2}
\dfrac{1 } {(q )_{\lambda_i }} \right\} 
=\sum_{\lambda \in \mathbb{Z}_{\geq}^n } \prod_{i=1}^{n} \left\{
\dfrac{ q^{m \lambda_i + \lambda_i^2 } }
{(q )_{\lambda_i }} \right\} \\
= \prod_{i=1}^n \left\{ \dfrac{ 1} {(q )_{\infty} } \right\} 
\cdot  \sum_{\lambda \in \mathbb{Z}^n }  
\prod_{i=1}^n \left\{ (-1)^{\lambda_i}\, q^{(-1/2+2m)\lambda_i +
   5\lambda_i^2/2 } \prod_{r=1 }^{m-1} (1-q^{r + \lambda_i}) \right\} 
\end{multline}
where the product $\prod_{r=1 }^{m-1} a_r=1$ for $m\leq 1$. Setting
 $m=\delta=0$ and $m=\delta=1$ and applying the \od\
 \Jtpi~(\ref{eq:JtpiBackground}) gives the result to be proved.
\end{proof}

Among all possible specializations $z_i =m/2+k(n-i)$
in~(\ref{eq:specializations}), the case $k=1$ gives similar
simplifications to the \od\ case when $m=0$ and $m=1$.  
Namely, the factors 
corresponding to the short and long roots in the $BC_n$
\RSi~(\ref{unspecializedRSi}) simplify in these cases. 

\begin{theorem}[$B_n$ and $D_n$ \RRis] 
\label{thm:DnRRis}
With the notation as above, $D_n$ multiple \RRis\ can be written as 
\begin{multline}
\label{DnRRis}
\sum_{\lambda \in \mathbb{Z}_{\geq}^n } \prod_{i=1}^{n} \left\{
\dfrac{ q^{(\delta+n-1) (\lambda_i-n+i) + (\lambda_i-n+i)^2 } }
{(q )_{\lambda_i }} \right\} \prod_{1\leq i<j \leq n}
   \left\{ (1-q^{\lambda_i-\lambda_j} )  \right\}  \\
= \dfrac{1}{2} (-1)^{\binom{n}{2}} \prod_{i=1}^n\left\{
   \dfrac{(q^5; q^5)_\infty}{(q)_\infty} q^{(n-i)(n-i+\delta/2)}
   \right\} \\ 
\cdot \det_{1\leq i,j\leq n} \bigg( q^{(j-1)(n-i+\delta/2)}\,
  \theta(q^{4n+2\delta+1-4i+j}; q^5 ) \bigg. 
   \bigg. \\ + q^{-(j - 1) (n-i+\delta/2)} \, \theta(q^{4n+2\delta
  +3-4i-j} ;q^5) \bigg)
\end{multline}
The cases $m=\delta=0$ and $m=\delta=1$ give the first and
the second $D_n$ \RRis, respectively. 

A single $B_n$ multiple \RRi\ which corresponds to $m=2$ may be
written, similarly, in the form  
\begin{multline}
\label{BnRRis}
\sum_{\lambda \in \mathbb{Z}_{\geq}^n } \prod_{i=1}^{n} \left\{
\dfrac{ q^{(1+n) (\lambda_i-n+i) + (\lambda_i-n+i)^2 } }
{(q)_{\lambda_i }} \right\} \prod_{1\leq i<j \leq n}
   \left\{ (1-q^{\lambda_i-\lambda_j} )  \right\}   \\
= (-1)^{\binom{n}{2}+n} \prod_{i=1}^n\left\{
   \dfrac{(q^5; q^5)_\infty}{(q)_\infty} q^{ (n-i+1/2)(n-i+1)
   }\right\} \\
\cdot \det_{1\leq i,j\leq n} \bigg( q^{(j-1/2)(n-i+1)} 
\theta(q^{6+4n-4i+j}; q^5 ) \\
- q^{-(j - 1/2)(n-i+1)} 
\theta(q^{7+4n-4i-j}; q^5 ) \bigg)  
\end{multline}
In both cases, $n$ is a positive integer and $\,\abs{q}<1$ as usual.  
\end{theorem}

\begin{proof}  
The $k=1$ case of~(\ref{specializedRSi_bothsides}) reduces to
\begin{multline}
\label{step1}
\sum_{\lambda \in \mathbb{Z}_{\geq}^n } \prod_{i=1}^{n} \left\{
\dfrac{ q^{(m+n-1) (\lambda_i-n+i) + (\lambda_i-n+i)^2 } }
{(q)_{\lambda_i }} \right\} 
\prod_{1\leq i<j \leq n}
   \left\{ (1-q^{\lambda_i-\lambda_j} )  \right\}   \\
= \prod_{i=1}^n \dfrac{1}{(q)_{\infty} } 
 \sum_{\lambda \in \mathbb{Z}^n }  
\prod_{i=1}^n \left\{ (-1)^{\lambda_i}\, q^{(-1/2+2m+3(n-i))\lambda_i +
   5\lambda_i^2/2 } \prod_{r=1 }^{m-1}  
(1-q^{r+n-i +\lambda_i})  \right\} \\ \cdot
\prod_{1\leq i< j \leq n} 
\left\{ (1-q^{-i+j + \lambda_i - \lambda_j}) 
(1- q^{ m + 2n-i-j +\lambda_i+\lambda_j} )
\right\} 
\end{multline}
where the product $\prod_{r=1 }^{m-1} a_r=1$ for $m\leq 1$. 

Setting $m=\delta\in \{0, 1\}$ and using the well--known determinant
evaluations~\cite{Krattenhaler1}
\begin{multline}
   \prod_{1\leq i<j \leq n} (1-x_i x_j^{-1}) (1-x_i x_j)  
   \\= \dfrac{1}{2} \, (-1)^{\binom{n}{2}} \prod_{i=1}^n x_i^{n-i}  
   \det_{1\leq i,j\leq n} \left( x_i^{j-1} + x_i^{-(j-1)} \right) 
\end{multline} 
the series on the \wps\ of~(\ref{step1}) becomes
\begin{multline}
  \dfrac{1}{2} \, (-1)^{\binom{n}{2}} \prod_{i=1}^n \left\{
  q^{(n-i)(n-i+\delta/2)} \right\}   
\det_{1\leq i,j\leq n} \bigg( \sum_{\lambda_i \in \mathbb{Z}}
   (-1)^{\lambda_i} q^{2(1+\delta +2n-2i) \lambda_i +5
     \binom{\lambda_i}{2} } 
   \bigg. \\ \bigg.  \cdot \left( q^{(j-1)(\lambda_i+n-i+\delta/2)}  +
     q^{-(j  - 1)(\lambda_i+n-i+\delta/2)} \right)  \bigg)
\end{multline}
The \od\ \Jtpi~(\ref{eq:JtpiBackground}) is then used to get the
desired $D_n$ identities.

Setting $m=2$ in~(\ref{step1}), and pursuing a similar line of thought 
using another determinant evaluation~\cite{Krattenhaler1} 
\begin{multline}
\prod_{1\leq i<j\leq n} (1-x_i x_j^{-1}) (1 - x_i x_j) \,
  \prod_{i=1}^n (1-x_i) \\ = (-1)^{\binom{n}{2}+n} \prod_{i=1}^n
  x_i^{n-i+1/2} \det_{1\leq 
  i, j \leq n} \left( x_i^{j-1/2} - x_i^{-(j-1/2)} \right) 
\end{multline}
gives the desired $B_n$ \RRi.
\end{proof} 

\begin{remark}
\label{DetFormRRis}
The Vandermonde determinant 
\begin{equation}
\label{Vander}
 \prod_{1\leq i<j\leq n} (1 - x_i x_j^{-1}) = (-1)^{\binom{n}{2}}
  \prod_{i=1}^n x_i^{1-i} \det_{1\leq i, j \leq n} \left( x_i^{n-j}
  \right)  
\end{equation}
may be used similarly to write the
balanced sides of~(\ref{step1}) as a determinant of theta functions as
well, which produces a determinant transformation identity where an 
$A_n$ type determinant equals a $B_n$ or $D_n$ type
determinant of theta functions. 

Using~(\ref{Vander}) and the recent result~(\ref{eq:GISgen}),
the determinant transformation identities
corresponding to~(\ref{DnRRis}) and~(\ref{BnRRis}) can be written in
the form 
\begin{multline}
\label{RRDetTransfIdenDn}
2 q^{\binom{n}{2}(1-n-3\delta/2) }
 \det_{1\leq i, j \leq n} \left( \pi_{\delta+i-j}(q) \right) \\
= \det_{1\leq i,j\leq n} \bigg( q^{(j-1)(n-i+\delta/2)}\,
  \theta(q^{4n+2\delta+1-4i+j}; q^5 ) \bigg. 
   \bigg. \\ + q^{-(j - 1) (n-i+\delta/2)} \, \theta(q^{4n+2\delta
  +3-4i-j} ;q^5) \bigg) 
\end{multline}
and 
\begin{multline}
\label{RRDetTransfIdenBn}
(-1)^{\binom{n}{2}} q^{-n(2n^2+3n-3)/4}
 \det_{1\leq i, j \leq n} \left( \pi_{2+i-j}(q) \right) \\
=  \det_{1\leq i,j\leq n} \bigg( q^{(j-1/2)(n-i+1)} 
\theta(q^{6+4n-4i+j}; q^5 ) \\
- q^{-(j - 1/2)(n-i+1)} 
\theta(q^{7+4n-4i-j}; q^5 ) \bigg)  
\end{multline}
respectively, where
\begin{equation}
\label{def:pi}
\pi_k(q):=
(-1)^k q^{-\binom{k}{2} }\theta(q^2; q^5) E_{k-2}(q) - (-1)^k
  q^{-\binom{k}{2} } \theta(q; q^5) D_{k-2}(q) 
\end{equation}
for $k\in\mathbb{Z}$, and the Schur polynomials $D_k(q)$ and $E_k(q)$
are defined as in~(\ref{eq:SchurPolyRecursive}). 

These results do not yield new relations between theta functions, yet
they appear to be new determinant transformation identities. A direct
proof using manipulations of theta functions and properties of
determinants is not obvious beyond dimension $n=2$. This is because
the entries of the Toeplitz matrix on the \lhs\ in all three cases
involve Schur polynomials of higher 
degrees as $n$ gets larger. It should be also noted that the
identities can be put into different forms by transposing matrices,
etc. 
\end{remark}

\begin{proof}
Using~(\ref{Vander}), the \rhs\ of the $k=1$ case
of~(\ref{specializedRSi_bothsides}) can be written in the form 
\begin{equation}
  (-1)^{\binom{n}{2}} \prod_{i=1}^{n} \left\{ 
    q^{(m-1+i)  (-n+i) }  \right\}
  \det_{1\leq i, j \leq n} \left( \sum_{\lambda_i \in \mathbb{Z}_{\geq}}
    \dfrac{ q^{(m+i-j)\lambda_i + \lambda_i^2 } }
    {(q )_{\lambda_i }} \right) 
\end{equation}
The identity~(\ref{eq:GISgen}) together with the
definition~(\ref{def:pi}) now imply that the \rhs\ equals 
\begin{equation}
(-1)^{\binom{n}{2}} \prod_{i=1}^{n} \left\{ 
q^{(m-1+i)  (-n+i) }  \right\}
 \det_{1\leq i, j \leq n} \left( \pi_{m+i-j}(q) \right) 
\end{equation}
The determinant transformation identities to be proved now follow.
\end{proof}

Multiple generalizations of other important \od\ $q$--series
identities can proved using the $BC_n$ \BL~\ref{OneParaBaileyLem} and
the Lemma~\ref{multilateralization}. In certain cases, it is possible to 
write such multiple generalizations as determinant transformation
identities involving theta functions similar to the ones given above
in~(\ref{RRDetTransfIdenDn}) and~(\ref{RRDetTransfIdenBn}). In other
cases, it is possible to find non--trivial product representations
producing genuine extensions of the classical results. The next
section discusses a remarkable such generalization of 
\Epnt. 

\subsection{Euler's Pentagonal Number Theorem}
An infinite family of $D_n$ \Epnt s can be written as follows.  
\begin{lemma}
Let $n$ be a positive and $k$ be a non--negative integer. Then
\begin{multline}
\label{Epnt}
(q)^n_\infty \!\!\! \prod_{1\leq i<j\leq n} \dfrac{(q^{k(j-i)})_\infty
  }{(q^{k(j-i+1)} )_{\infty} } 
= \sum_{\mu\in\mathbb{Z}^n } 
 (-1)^{|\mu|} q^{-kn(\mu)+3n(\mu')+|\mu|(k(n-1)+1)} \\ 
\cdot \sideset{}{^{k,1}}\prod_{1\leq i<j \leq n}
 \left\{ \dfrac{(q^{1+k(j-i-1) + \mu_i-\mu_j})_{\infty} }  
{(q^{k(1+2n-i-j) + \mu_i+\mu_j})_{\infty} } 
 \dfrac{ (q^{1+k(-1+2n-i-j) + \mu_i+\mu_j})_{\infty} }
{(q^{k(j-i+1) + \mu_i -\mu_j})_{\infty} } \right\}
\end{multline}
where $q\in\mathbb{C}$ such that $\abs{q}<1$ as usual. 
\end{lemma}

\begin{proof}
\label{proofEpnt}
For a rectangular partition $\lambda=k^n$, the terminating \sfs\
summation~(\ref{sfs}) can be written in the form 
\begin{multline}
\dfrac{ (qb, qb/\rho_1 \sigma_1)_{k^n} } {(qb/\sigma_1,
  qb/\rho_1)_{k^n} } = \sum_{\mu\subseteq \lambda}  
  t^{2n(\mu)+(1-n)|\mu|}  \left( \dfrac{q^{1+k}b}{\sigma_1\rho_1}
  \right)^{|\mu|} \\ \cdot \prod_{i=1}^n
  \left\{\dfrac{(1-bt^{2-2i}q^{2\mu_i})} {(1-b t^{2-2i})}
    \right\} 
  \dfrac{(bt^{1-n}, \sigma_1, \rho_1, q^{-k})_{\mu}}{(qt^{n-1},
  qb/\sigma_1, qb/\rho_1, q^{1+k}b)_{\mu}} \\
  \cdot \prod_{1\leq i<j \leq  n} \left\{\dfrac{ (qt^{j-i})_{\mu_i-\mu_j}
  (bt^{3-i-j})_{\mu_i+\mu_j}} {(qt^{j-i-1})_{\mu_i-\mu_j} (b
  t^{2-i-j})_{\mu_i+\mu_j} } \dfrac{(t^{j-i+1})_{\mu_i -\mu_j}
  (qbt^{2-i-j})_{\mu_i+\mu_j}} {(t^{j-i})_{\mu_i -\mu_j}
  (qbt^{1-i-j})_{\mu_i+\mu_j}} \right\}
\end{multline}
by the virtue of the \Wdegf~(\ref{conj:degform}). Sending $\rho_1$,
$\sigma_1$ and $k$ to $\infty$ 
using the Dominated Convergence Theorem, similar to the
proof of the \RSi\ in Lemma~\ref{lem:unspecializedRSi}, gives
\begin{multline}
\label{Epnt1b}
(qb)_{\infty^n} = \sum_{\mu, \ell(\mu)\leq n} 
  t^{-n(\mu)+(1-n)|\mu|} b^{|\mu|} (-1)^{|\mu|} q^{|\mu|+3n(\mu')} \\ 
  \cdot \prod_{i=1}^n  
  \left\{\dfrac{(1-bt^{2-2i}q^{2\mu_i})} {(1-b t^{2-2i})}
    \right\} \dfrac{(bt^{1-n})_\mu }{ (qt^{n-1})_\mu } \\
\cdot \prod_{1\leq i<j \leq  n} \left\{\dfrac{ (qt^{j-i})_{\mu_i-\mu_j}
  (bt^{3-i-j})_{\mu_i+\mu_j}} {(qt^{j-i-1})_{\mu_i-\mu_j} (b
  t^{2-i-j})_{\mu_i+\mu_j} } \dfrac{(t^{j-i+1})_{\mu_i -\mu_j}
  (qbt^{2-i-j})_{\mu_i+\mu_j}} {(t^{j-i})_{\mu_i -\mu_j}
  (qbt^{1-i-j})_{\mu_i+\mu_j}} \right\}
\end{multline}
Now let $t=q^k$ and $b = q^{m+2k(n-1)}$. 
This corresponds to setting $q^{z_i} = b^{1/2} t^{1-i}$ for $i\in [n]$
and specializing $z_i = m/2+k(n-i)$ for  non--negative integers $m$
and $k$ as in~(\ref{eq:specializations}). Flip appropriate terms using
the definition of the \qPs~(\ref{qPochSymbol}) and simplify to get  
\begin{multline}
(q)^n_\infty \!\!\! \prod_{1\leq i<j\leq n} \!\!
  \dfrac{(q^{k(j-i)})_\infty }{(q^{k(j-i+1)} )_{\infty} } 
= \sum_{\mu, \ell(\mu)\leq n} 
 (-1)^{|\mu|} q^{-kn(\mu)+3n(\mu')+|\mu|(m+k(n-1)+1)} \\  \cdot
 \prod_{i=1}^n 
\left\{\dfrac{(q^{m+2k(n-i) +2\mu_i})_{\infty} }
 {(q^{1+m+2k(n-i)+2\mu_i})_{\infty} }  
\dfrac{ (q^{1+k(n-i) +\mu_i})_{\infty} } {(q^{m + k(n-i) +\mu_i}
 )_{\infty}}  \right\} \\ 
\cdot \sideset{}{^{k,\aleph_\lambda}}\prod_{1\leq i<j \leq n}
 \left\{\dfrac{(q^{k(j-i) + \mu_i -\mu_j})_{\infty} }  
{ (q^{1+k(j-i) + \mu_i-\mu_j} )_{\infty}} 
\dfrac { (q^{m+k(2n-i-j) + \mu_i+\mu_j})_{\infty} } 
{ (q^{1+m+k(2n-i-j) + \mu_i+\mu_j})_{\infty} } \right. \\ \left. 
\cdot \dfrac{(q^{1+k(j-i-1) + \mu_i-\mu_j})_{\infty} }  
{(q^{m+k(1+2n-i-j) + \mu_i+\mu_j})_{\infty} } 
 \dfrac{ (q^{1+m+k(-1+2n-i-j) + \mu_i+\mu_j})_{\infty} }
{(q^{k(j-i+1) + \mu_i -\mu_j})_{\infty} } \right\}
\end{multline}
With this specialization, Remark~\ref{HyperoctahedralSymmetry} shows
that all the conditions of the multilateralization
lemma~\ref{multilateralization} are satisfied. Applying this lemma,
therefore, yields the identity 
\begin{multline}
(q)^n_\infty \!\!\! \prod_{1\leq i<j\leq n} \!\!
  \dfrac{(q^{k(j-i)})_\infty }{(q^{k(j-i+1)} )_{\infty} } 
= \sum_{\mu\in\mathbb{Z}^n }  
 (-1)^{|\mu|} q^{-kn(\mu)+3n(\mu')+|\mu|(m+k(n-1)+1)} \\  \cdot
\sideset{}{^{m,1}} \prod_{i=1}^n 
\left\{ \dfrac{ (q^{1+k(n-i) +\mu_i})_{\infty} } {(q^{m + k(n-i)
      +\mu_i} )_{\infty}}  \right\} \\ 
\cdot \sideset{}{^{k,1}}\prod_{1\leq i<j \leq n}
 \left\{ \dfrac{(q^{1+k(j-i-1) + \mu_i-\mu_j})_{\infty} }  
{(q^{m+k(1+2n-i-j) + \mu_i+\mu_j})_{\infty} } 
 \dfrac{ (q^{1+m+k(-1+2n-i-j) + \mu_i+\mu_j})_{\infty} }
{(q^{k(j-i+1) + \mu_i -\mu_j})_{\infty} } \right\}
\end{multline}
Finally, mimicing the proof of the classical \Epnt\ and setting $m=0$
(i.e., $t=q^k$ and $b = q^{2k(n-1)}$) gives the desired identity. 
Notice that the resulting identity~(\ref{Epnt}) could also be written
in the form 
\begin{multline}
\prod_{i=1}^{n} \left\{ (q)_{\infty} (1-q^{ki})^{n-ki} \right\} 
= \sum_{\mu\in\mathbb{Z}^n } 
 (-1)^{|\mu|} q^{-kn(\mu)+3n(\mu')+|\mu|(k(n-1)+1)} \\ 
\cdot \prod_{1\leq i<j \leq n}
 \left\{ \prod_{r=0}^{2k-2} (1 - q^{r+1+k(-1+j-i) +
       \mu_i-\mu_j}) (1-q^{r+1+k(-1+2n-i-j) +
       \mu_i+\mu_j}) \right\}
\end{multline}
using the standard product notation instead.
\end{proof}

\begin{remark}
First note that, as in Remark~\ref{centerz}, the result may be extended
into negative integers $m$ and $k$ with some extra work.

Similar to the alternative forms~(\ref{RRDetTransfIdenDn})
and~(\ref{RRDetTransfIdenBn}) of the $D_n$ and $B_n$ \RRis, $D_n$
\Epnt~(\ref{Epnt}) can be written as a determinant identity, in
particular, in the case $k=1$.
   
Applying same steps in the proof of Remark~\ref{DetFormRRis}
to the $k=1$ instance of~(\ref{Epnt}) gives the identity
\begin{multline}
\label{Epnt_k1}
2 \, (-1)^{\binom{n}{2}} q^{-n(n+1)(2n+1)/6} \,\, \theta(q; q^3)^n
\!\!\! \prod_{1\leq i<j\leq n} \dfrac{(q^{j-i})_\infty }{(q^{j-i+1}
  )_{\infty} }  \\ 
= \det_{1\leq i,j\leq n} \left( q^{(n-i) (j-1)} \theta(q^{(2n-2i+j)}; q^3) 
   +  q^{-(n-i) (j-1)} \theta(q^{(2n-2i-j+2)}; q^3)  \right) 
\end{multline}
The identity~(\ref{Epnt_k1}) does not, again, give rise to new relations
between theta functions. Yet, it appears to be a new determinant
identity involving theta functions. It should be also noted that the
identity can be put into different forms by transposing matrices, etc.
\end{remark}

\Epnt\ and \RRis\ correspond to the first and second iteration of the
\BL. In general, iterating Bailey Lemma $N$ times gives rise to
so--called generalized \Wt\ and a limiting case of it yields the
generalized \RSi. The \Jtpi\ is then used to compute a product
representation under certain specializations in the classical
case. The resulting identities are called the (extreme cases of)
\AGis\ which generalize the \RRis. The full \AGis\ may be proved using
the classical two parameter \BL~\cite{AgarwalA1}.

A multiple analogue~(\ref{generalWatson}) of the generalized \Wt\ is
already obtained above by iterating the $BC_n$ Bailey
Lemma~\ref{OneParaBaileyLem}. A limiting case of this transformation
gives a multiple extension of the extreme cases of \AGis\ similar to
the \RRis\ proved in this paper. Both the extreme cases and a
remarkable full version of \AGis\ are investigated 
in another publication~\cite{Coskun4}.         

The $BC_n$ \BL s~\ref{OneParaBaileyLem} and~\ref{twoParameterBL} may
be used to prove generalizations of many important 
identities as in the classical case. Elliptic and basic
(trigonometric) root system analogues of several significant classical
summation and transformation identities are proved using these results
in~\cite{Coskun3} (also see~\cite{Rosengren1}, \cite{RosengrenS1},
\cite{Warnaar2}, \cite{Schlosser1}).  The classical one parameter \BL\
has another extension known as WP--\BL~\cite{Andrews1}, \cite{AndrewsB1},
\cite{Spiridonov1}. An interpretation of the elliptic $BC_n$
\BL~\ref{twoParameterBL} in the 
setting of a multivariate interpolation problem  
is used in~\cite{Coskun5} to obtain what is
called interpolation $BC_n$ \BL, which may be considered as a
multiple elliptic generalization of the WP--\BL.

\end{document}